\input amstex
\magnification=1200
\loadmsam
\loadmsbm
\loadeufm
\loadeusm
\UseAMSsymbols
\hsize=6.00 true in
\vsize=8.7 true in

\input pictex.tex

\font\thinlinefont=cmr5
\font\BIGtitle=cmss10 scaled\magstep2
\font\authorfont=cmss10 scaled\magstep1
\font\proclaimfont=cmssdc10
\font\sectionfont=cmss10 scaled\magstep1
\font\sans=cmss10


\def\scr#1{{\fam\eusmfam\relax#1}}

\def\scrC{{\scr C}}
\def\scrS{{\scr S}}
\def\scrW{{\scr W}}

\def\bfone{\bold{1}}
\def\bftwo{\bold{2}}
\def\bfA{\bold{A}}
\def\bfK{\bold{K}}
\def\bfL{\bold{L}}
\def\bfb{\bold{b}}
\def\bfp{\bold{p}}
\def\bfq{\bold{q}}
\def\Ch{\bold{Ch}}

\def\Ftil{\widetilde{F}}

\def\Ho{\text{\rm Ho}}
\def\hobaw{\text{\rm holim}^{BW}}
\def\hofoo{\text{\rm holim}^{T}}
\def\hocofoo{\text{\rm hocolim}^{T}}
\def\Sd{\text{\rm Sd}}
\def\we{\text{\rm we}}
\def\fib{\text{\rm fib}}
\def\op{\text{\rm op}}
\def\id{\text{\rm id}}
\def\Tot{\text{\rm Tot}}
\def\Map{\text{\rm Map}}
\def\ho{\text{\rm ho}}
\def\hoend{\text{\rm ho}\!\int}
\def\holim{\operatornamewithlimits{%
   \text{\rm holim\,}}}
\def\Cat{\text{\rm Cat}}
\def\HOM{\underline{\Hom}}
\def\Hom{\text{\rm Hom}}
\def\hom{\text{\rm hom}}
\def\Sets{\text{$\bold{Sets}$}}
\def\fib{\text{\rm fib}}
\def\cof{\text{\rm cof}}

\def\nspace{\lineskip=1pt\baselineskip=12pt%
     \lineskiplimit=0pt}
\def\dspace{\lineskip=2pt\baselineskip=18pt%
     \lineskiplimit=0pt}
\def\Proclaim#1{\medbreak\noindent%
     {\proclaimfont #1.\enspace}\it\ignorespaces}
\def\finishproclaim{\par\rm
     \ifdim\lastskip<\medskipamount\removelastskip
     \penalty55\medskip\fi}
\def\proof{\par\noindent {\it Proof:}\enspace}
\def\proofof#1{\par\medskip\noindent{\it Proof}
     \enspace}
\def\references#1{\par
  \centerline{\sectionfont References}\bigskip
     \parindent=#1pt\nspace}
\def\Ref[#1]{\par\medskip\hang\indent%
     \llap{\hbox to\parindent
     {[#1]\hfil\enspace}}\ignorespaces}
\def\Item#1{\par\smallskip\hang\indent%
     \llap{\hbox to\parindent
     {#1\hfill$\,\,$}}\ignorespaces}
\def\ItemItem#1{\par\indent\hangindent2%
     \parindent \hbox to\parindent%
     {#1\hfill\enspace}\ignorespaces}
\def\updot{{\dsize\cdot}}

\def\Le{{\mathchoice{\,{\scriptstyle\le}\,}
  {\,{\scriptstyle\le}\,}%
  {\,{\scriptscriptstyle\le}\,}%
  {\,{\scriptscriptstyle\le}\,}}}
\def\Ge{{\mathchoice{\,{\scriptstyle\ge}\,}
  {\,{\scriptstyle\ge}\,}
  {\,{\scriptscriptstyle\ge}\,}%
  {\,{\scriptscriptstyle\ge}\,}}}

\def\arrowsim{\,\smash{\mathop{\longrightarrow}%
\limits^{\lower1.5pt\hbox{$\scriptstyle\sim$}}}\,}

\def\:{\colon\,}

\def\maprightarrowtail#1{\,\,\raise7pt\hbox{%
  $\mathop{\vtop{\ialign{##\crcr
\noalign{\nointerlineskip}\crcr
  \hfil\lower2pt\hbox{${\sssize #1}$}\hfil\crcr
\noalign{\nointerlineskip}\crcr
  \hfil\raise2pt\hbox{$\rightarrowtail$}\crcr
  \noalign{\nointerlineskip}\crcr}}}$}\,\,}

\def\longmaprightarrowtail#1{\,\,\raise7pt\hbox{%
  $\mathop{\vtop{\ialign{##\crcr
\noalign{\nointerlineskip}\crcr
  \hfil\lower2pt\hbox{${\sssize #1}$}\hfil\crcr
\noalign{\nointerlineskip}\crcr
  \hfil\raise2pt\hbox{$\longrightarrowtail$}\crcr
  \noalign{\nointerlineskip}\crcr}}}$}\,\,}

\def\limright{\mathop{\vtop{\ialign{##\crcr
    \hfil\rm lim\hfil\crcr
    \noalign{\nointerlineskip}
\lower5.5pt\hbox to 15pt{\rightarrowfill}\crcr
    \noalign{\nointerlineskip}\crcr}}}}

\def\limleft{\mathop{\vtop{\ialign{##\crcr
  \hfil\rm lim\hfil\crcr
  \noalign{\nointerlineskip}
\lower5.5pt\hbox to 15pt{\leftarrowfill}\crcr
  \noalign{\nointerlineskip}\crcr}}}}

\font\eightrm=cmr8
\font\sevenrm=cmr7
\font\sixrm=cmr6
\font\fiverm=cmr5
\font\eighti=cmmi8
\font\sixi=cmmi6
\font\fivei=cmmi5
\font\eightsy=cmsy8
\font\sixsy=cmsy6
\font\fivesy=cmsy5
\font\tenex=cmex10
\font\eightit=cmti8
\font\eightsl=cmsl8
\font\eighttt=cmtt8
\font\eightbf=cmbx8
\font\sixbf=cmbx6
\font\fivebf=cmbx5

\def\eightpoint{\def\rm{\fam0\eightrm}
  \textfont0=\eightrm \scriptfont0=\sixrm \scriptscriptfont0=\fiverm
  \textfont1=\eighti \scriptfont1=\sixi \scriptscriptfont1=\fivei
  \textfont2=\eightsy \scriptfont2=\sixsy \scriptscriptfont2=\fivesy
  \textfont3=\tenex \scriptfont3=\tenex \scriptscriptfont3=\tenex
  \textfont\itfam=\eightit  \def\it{\fam\itfam\eightit}%
  \textfont\slfam=\eightsl  \def\sl{\fam\slfam\eightsl}%
  \textfont\ttfam=\eighttt  \def\tt{\fam\ttfam\eighttt}%
  \textfont\bffam=\eightbf  \scriptfont\bffam=\sixbf
   \scriptscriptfont\bffam=\fivebf  \def\bf{\fam\bffam\eightbf}%
  \normalbaselineskip=9pt
  \setbox\strutbox=\hbox{\vrule height7pt depth2pt width0pt}%
  \let\sc=\sixrm  \normalbaselines\rm}


\catcode`\@=11 
\def\footmarkform@#1{$\m@th^{#1}$} 
\let\thefootnotemark\footmarkform@ 
\def\makefootnote@#1#2{\insert\footins 
{\interlinepenalty\interfootnotelinepenalty 
\eightpoint\splittopskip\ht\strutbox\splitmaxdepth\dp\strutbox 
\floatingpenalty\@MM\leftskip\z@\rightskip\z@%
     \spaceskip\z@\xspaceskip\z@
 \leavevmode{#1}\footstrut\ignorespaces#2\unskip\lower\dp\strutbox
 \vbox to\dp\strutbox{}}}
\newcount\footmarkcount@
\footmarkcount@\z@
\def\footnotemark{\let\@sf\empty\relaxnext@
 \ifhmode\edef\@sf{\spacefactor\the\spacefactor}\/\fi
 \DN@{\ifx[\next\let\next@\nextii@\else
  \ifx"\next\let\next@\nextiii@\else
  \let\next@\nextiv@\fi\fi\next@}%
 \DNii@[##1]{\footmarkform@{##1}\@sf}%
 \def\nextiii@"##1"{{##1}\@sf}%
 \def\nextiv@{\iffirstchoice@\global\advance\footmarkcount@\@ne\fi
  \footmarkform@{\number\footmarkcount@}\@sf}%
 \FN@\next@}
\def\footnotetext{\relaxnext@
 \DN@{\ifx[\next\let\next@\nextii@\else
  \ifx"\next\let\next@\nextiii@\else
  \let\next@\nextiv@\fi\fi\next@}%
 \DNii@[##1]##2{\makefootnote@{\footmarkform@{##1}}{##2}}%
 \def\nextiii@"##1"##2{\makefootnote@{##1}{##2}}%
 \def\nextiv@##1{\makefootnote@{\footmarkform@{\number%
     \footmarkcount@}}{##1}}%
 \FN@\next@}
\def\footnote{\let\@sf\empty\relaxnext@
 \ifhmode\edef\@sf{\spacefactor\the\spacefactor}\/\fi
 \DN@{\ifx[\next\let\next@\nextii@\else
  \ifx"\next\let\next@\nextiii@\else
  \let\next@\nextiv@\fi\fi\next@}%
 \DNii@[##1]##2{\footnotemark[##1]\footnotetext[##1]{##2}}%
 \def\nextiii@"##1"##2{\footnotemark"##1"\footnotetext"##1"{##2}}%
 \def\nextiv@##1{\footnotemark\footnotetext{##1}}%
 \FN@\next@}
\def\adjustfootnotemark#1{\advance\footmarkcount@#1\relax}
\def\footnoterule{\kern-3\p@
  \hrule width 5pc\kern 2.6\p@} 
\catcode`\@=\active

\def\lmapdown#1{\Big\downarrow\llap%
   {$\vcenter{\hbox{$\scriptstyle#1\enspace$}}$}}

\def\rmapdown#1{\Big\downarrow\kern-1.0pt%
  \vcenter{\hbox{$\scriptstyle#1$}}}

\def\mapright#1{\smash{\mathop{\longrightarrow}%
     \limits^{#1}}}

\def\longmapright#1#2{\smash{\mathop
    {\hbox to #1pt{\rightarrowfill}}\limits^{#2}}}

\def\longmaprightsub#1#2{\mathop
{\hbox to #1pt{\rightarrowfill}}\limits_{#2}}

\def\longmaprightsubsup#1#2#3{%
  \mathop{\hbox to #1pt{\rightarrowfill}}%
 \limits_{\raise3pt\hbox{$\scriptstyle #2$}}^{#3}}

\def\longmapleft#1#2{\smash{\mathop
    {\hbox to #1pt{\leftarrowfill}}\limits^{#2}}}

\def\Rmapdown#1{\Bigg\downarrow%
  \rlap{$\kern-1.0pt\vcenter{%
  \hbox{$\scriptstyle#1$}}$}}

\def\Lmapdown#1{\rlap{$\vcenter{%
   \hbox{$\scriptstyle#1$}}$}\kern4.0pt
  \bigg\downarrow}

\def\longrightarrowtail{\,\raise2.50pt
  \hbox{$\sssize >$}\kern-2.0pt\raise1.2pt
  \hbox{$\longrightarrow$}\,}

\def\longleftarrowtail{\,\raise1.2pt
  \hbox{$\longleftarrow$}\kern-2.0pt
  \raise2.50pt\hbox{$\sssize <$}\,}

\def\twoheadSimArrow#1{\raise5.0pt\hbox{%
  $\mathop{\vtop{\ialign{##\crcr 
\noalign{\nointerlineskip}\crcr
  \hfil$\!\!\sssize\sim$\hfil\crcr
\noalign{\nointerlineskip}\crcr
 \hfil\null\raise5.29pt\hbox to #1pt{\hrulefill}%
\!\raise3.0pt\hbox{$\twoheadrightarrow$}\hfil\crcr
    \noalign{\nointerlineskip}\crcr}}}$}}

\def\maptwohead#1#2{\,\raise5.0pt\hbox{%
  $\mathop{\vtop{\ialign{##\crcr
\noalign{\nointerlineskip}\crcr
  \hfil$\!\!{\ssize #2}$\hfil\crcr
\noalign{\nointerlineskip}\crcr
  \hfil\null\raise5.29pt\hbox to #1pt{\hrulefill}%
\!\raise3.0pt\hbox{$\twoheadrightarrow$}\hfil\crcr
  \noalign{\nointerlineskip}\crcr}}}$}\,}

\def\lrmapdown#1#2{\llap{$\vcenter{%
  \hbox{$\scriptstyle#1$}}$}
 \Big\downarrow\rlap{$\vcenter{%
  \hbox{$\scriptstyle#2$}}$}}

\def\twoheaddownarrow{%
  \raise7pt\hbox{$\big\vert$}
\hbox{\null\kern-4.25pt%
  \lower2pt\hbox{$\downarrow$}%
 \lower4pt\hbox{\null\kern-5pt$\downarrow$}}}

\def\downarrowtail{\raise11.50pt%
  \hbox{$\sssize\vee$}\kern-5.75pt\Big\downarrow}

\def\uparrowtail{\Big\uparrow
  \kern-5.65pt\lower9.25pt\hbox{$\sssize\wedge$}}

\def\Rtwoheadmapdown#1{%
  \raise7pt\hbox{$\big\vert$}
  \hbox{\null\kern-4.24pt%
  \lower2pt\hbox{$\downarrow$}%
  \lower4pt\hbox{\null\kern-5.00pt$\downarrow$}}
  \rlap{$\kern-1.0pt
  \vcenter{\hbox{\raise5pt\hbox{$\ssize#1$}}}$}}

\def\Ldownarrowtail#1{\raise11.50pt
  \hbox{$\sssize\vee$}\kern-5.75pt\Big\downarrow
  \llap{$\vcenter{\hbox{$\ssize#1\enspace$}}$}}

\def\Rdownarrowtail#1{\raise11.50pt
\hbox{$\sssize\vee$}\kern-5.75pt\Big\downarrow
\rlap{$\kern-1.0pt\vcenter{\hbox{$\ssize#1$}}$}}

\parindent=25pt
\document
\footline={\hfil}

\centerline{\BIGtitle 
Homotopy Ends and Thomason Model Categories}

\bigskip
\line{\hfill\vtop{\hbox{\authorfont  
Charles Weibel\raise6pt\hbox{\sevenrm 1}}
\smallskip
\hbox{Institute for Advanced Study}
\hbox{Princeton, NJ \ 08540}
\hbox{weibel\@math.ias.edu}}\hfill}
\footnote""{\null\kern-25pt
$^1$Partially supported by NSF grants}

\bigskip\bigskip
\dspace
In the last year  of his life, Robert W.~Thomason reworked the notion
of a model category, used to adapt homotopy theory to algebra,
and used homotopy ends to affirmatively solve a problem raised by
Grothendieck: find a notion of model structure which is inherited by
functor categories.
The axioms for such a Thomason model category were published later in
[WT]. In this paper we explain and prove Thomason's results, based on
his private notebooks [T78]--[T86].

Each of Thomason's 180-page notebooks are
spiral-bound with an index on the inside cover.
We have chosen to follow their internal
citation method: our citation [T$x,y$]
refers to page $y$ of Thomason's notebook $x$.

We first present Thomason's ideas about homotopy ends and its
generalizations (sections 1 and 2). 
These are first formulated for complete simplicial closed model
categories in the sense of Quillen [QH], because of their
usefulness in this setting, as demonstrated by Dwyer and Kan [DK].

Thomason's axioms and examples come next,
followed by a proof that the homotopy
category $\ho\scrC$ of a Thomason model
category $\scrC$ exists (sections 3-5).
In the last two sections (6--7) we prove the main theorem: 
the functor category $\scrC^{\bfK}$ inherits a Thomason model
structure, at least when $\scrC$ is right
enriched over simplicial sets, and
fibrations are preserved by products and inverse limits.
The result was previously known when $\scrC$ was the category of
simplicial sets (Quillen, Joyal and Jardine) and also when $\bfK$ is any
partially ordered set (Grothendieck).

\newpage
\footline={\hss\tenrm \folio\hss}
\pageno=1

\centerline{\sectionfont \S1 Homotopy Ends}

\bigskip\noindent
Recall from [Mac] that the {\it end} of a
functor $F\colon\bfK^{\op}\times\bfK\to\scrC$
is an object $e=\int\nolimits_{\bfK}F(K,K)$
together with maps $e\to F(K,K)$ which are
compatible in the sense that for each
$K_1\to K_2$ the evident square commutes:
$$
\CD
e @>>> F(K_1,K_1)\\
@VVV @VVV\\
F(K_2,K_2) @>>> F(K_1,K_2)\,.
\endCD
$$
Moreover, $e$ is universal with respect to this property.
If $\bfK$ is a small category and $\scrC$ is complete, then 
ends always exist, and commute with limits [Mac, IX.5].

The ``Tot'' construction for cosimplicial
objects is an example of
an end in $\scrC$.
Its construction requires that for every
finite simplicial set $K$ there is a
{\it mapping object functor} $\Map(K,-)\colon\,
\scrC\to\scrC$ which is natural in $K$.
Given a cosimplicial object $X^{\updot}$
in $\scrC$, we have a functor
$\Delta^{\op}\times\Delta\to\scrC$ defined
by $(p,q)\mapsto \Map(\Delta[p],X^q)$.

\Proclaim{Definition 1.1} \rm
If $X^{\updot}$ is a cosimplicial object
in a complete category $\scrC$, and
$\scrC$ is equipped with a mapping object
functor, we define the {\it total object}
of $X^{\updot}$ to be the end
$$
\Tot(X^{\updot})=\int\nolimits_{\Delta}
\Map(\Delta[p],X^p).
$$
\finishproclaim

In case $\scrC$ is a simplicial closed
model category in the sense of Quillen
[QH], it is well known that this
definition of $\Tot$ agrees with the
Bousfield-Kan definition in [BK].

In a model category, it is not generally
true that an end will preserve weak
equivalences.
To remedy this, we introduce the notion of
a ``derived end.''

Given $F\colon\,\bfK^{\op}\times
\bfK\to\scrC$, with $\scrC$ complete and
$\bfK$ small, the {\it cosimplicial
replacement} of $F$ is the cosimplicial
object $\prod^*F$ of $\scrC$ defined by
$$
p\longmapsto \prod\nolimits^p F=\prod\limits_{
K_0\to\cdots\to K_p}F(K_0,K_p),
$$
where the indexing set runs over 
all $p$-tuples $K_0\to K_1\to\cdots\to K_p$ of
composable maps in $\bfK$, i.e.,
all functors $\bfp\to\bfK$.
This is equivalent to the Bousfield-Kan
simplicial replacement of $F$ [BK, XI.5],
except for a change in orientation due to
a different orientation of the nerve of
$\bfK$ on p.~291 [BK].

\Proclaim{Lemma 1.2}
The end $\int\limits_{\bfK}F(K,K)$ is the
equalizer 
of $~~\prod\nolimits^0 F\rightrightarrows
\prod\nolimits^1 F$.
\finishproclaim

\proof {\rm ([T78, 115])}
Let $\pi^0 F$ denote the equalizer.
The projections $\pi^0
F\to\prod\nolimits^0 F\to F(K,K)$ have the
property that for every map $K_0\to K_1$
in $\bfK$ the diagram
$$
\CD
\pi^0F @>>>  F(K_0,K_0)\\
@VVV @VVV\\
F(K_1,K_1) @>>> F(K_0,K_1)
\endCD
$$
commutes.
This gives a map from $\pi^0F$ to the end
$\int_{\bfK}F (K,K)$.
Since the end also equalizes
$\prod\nolimits^0F\rightrightarrows
\prod\nolimits^1 F$, this must be an
isomorphism.\qquad
$\square$

\Proclaim{Definition 1.3} \rm
If $\scrC$ is a complete category equipped
with a mapping object functor, the {\it
homotopy end} $\,\,\,\hoend F$ of
$F\colon\,\bfK^{\op}\times\bfK\to\scrC$ is
the total object of $\prod\nolimits^*F$:
$$
\hoend
F=\Tot\left(\prod\nolimits^* F\right)=
\int\nolimits_{\Delta}\Map\left(\Delta[p],
\prod\nolimits^pF\right)\,\,.
$$
\finishproclaim

\Proclaim{Compatibility 1.4} 
\rm Suppose that $F$ factors as the projection
$\bfK^{\op}\times\bfK\to\bfK^{\op}$
followed by a functor
$\Ftil\colon\,\bfK^{\op}\to\scrC$.
Then the homotopy end of $F$ is just the
Bousfield-Kan homotopy limit of $\Ftil$ as
described in 4.5 and 5.2 of [BK, XI]:
$$
\hoend F=\holim \Ftil\,\,.
$$
Since we have used $\bfK^{\op}$, there is
no orientation problem.

\Proclaim{Functoriality 1.5} \rm
A natural transformation
$\eta\colon\,F\Rightarrow G$  of functors
$\bfK^{\op}\times\bfK\to\scrC$ induces
a map
$\prod\nolimits^*F\to\prod\nolimits^*G$,
and so a natural map $\,\,\,\hoend F\to\hoend G$.
The homotopy end functor
$\Cat(\bfK^{\op}\times\bfK,\scrC)\to\scrC$
commutes with all limits (because $\Tot$
and $F\mapsto \prod\nolimits^* F$ preserve
limits).
Given a functor $\Phi\colon\,\bfL\to\bfK$
and $F\colon\,\bfK^{\op}\times
\bfK\to\scrC$, the diagram
$$
\CD
\Cat(\bfK^{\op}\times\bfK,\scrC)
@>{\Phi}>>
\Cat(\bfL^{\op}\times\bfL,\scrC)\\
\null\kern40pt\raise10pt\hbox{$\ssize\hoend$}\kern-12pt
\vbox{\beginpicture
\setcoordinatesystem units <.50cm,.50cm>
\linethickness=1pt
\setshadesymbol ({\thinlinefont .})
\setlinear
%
%
\linethickness= 0.500pt
\setplotsymbol ({\thinlinefont .})
\plot  6.032 20.955  8.572 18.733 /
%
%
\plot  8.340 18.852  8.572 18.733  8.423 18.948 /
\linethickness=0pt
\putrectangle corners at  6.007 20.980 and  8.598 18.707
\endpicture}
&&\null\kern-40pt\vbox{\beginpicture
\setcoordinatesystem units <.50cm,.50cm>
\linethickness=1pt
\setshadesymbol ({\thinlinefont .})
\setlinear
%
%
\linethickness= 0.500pt
\setplotsymbol ({\thinlinefont .})
\plot 10.795 21.273  8.572 18.733 /
%
%
\plot  8.692 18.965  8.572 18.733  8.788 18.882 /
\linethickness=0pt
\putrectangle corners at  8.547 21.298 and 10.820 18.707
\endpicture}\null\kern-15pt
\raise10pt\hbox{$\ssize\hoend$}\\
\vspace{-5pt}
&\scrC&
\endCD
$$
commutes up to a natural $2$-cell
$\varphi$ which is induced by the
``projection''
$\prod\nolimits_{\bfK}^*F\to
\prod\nolimits_{\bfL}^*F$: the component
indexed by $L\: p\to\bfL$ comes from the
component indexed by $\Phi L$.

We now introduce homotopy theory into the
discussion.
Suppose to fix ideas that $\scrC$ is a
simplicial closed model category in the
sense of Quillen [QH, II.2.2], and that
$\scrC$ is complete.
We first observe that $\scrC$ is
cotensored over all simplicial sets, not
just finite ones, as in Quillen's axiom (SM0).

\Proclaim{Lemma 1.6 {\rm [T78, 58]}}
There is a 
{\rm mapping object functor}
$\Map(K,-)\colon\,\scrC\to\scrC$ for every
simplicial set $K$, which is
natural in $K$, and an isomorphism
$$
\Hom_{\scrC}(X,\Map(K,Y))\cong
\Hom_{\Delta^{\op}\bold{Sets}}
(K,\hom(X,Y)),
$$
which is natural in $K$, $X$ and $Y$.

Quillen's axiom (SM7) holds in this
setting.
If $i\colon\, K\to L$ is a cofibration in
$\Delta^{\op}\Sets$ (i.e., an injection)
and $f\colon\,X\to Y$ is a fibration in
$\scrC$, then
$$
\Map(L,X)\to\Map(K,X)\times_{\Map(K,Y)}
\Map(L,Y)
$$
is a fibration in $\scrC$, and is a weak
equivalence if either $i$ or $f$ is.
Moreover, the functor
$\Map(\,\,\,,\,\,\,)\: 
\Delta^{\op}\Sets\times\scrC\to\scrC$
preserves limits in $\scrC$ and converts
colimits of simplicial sets to limits in
$\scrC$.
\finishproclaim

\proof
Write $K=\sqcup K_\alpha$, where $K_\alpha$
runs over all finite subspaces of $K$, and
set $\Map(K,X)=\lim\,\Map(K_\alpha,X)$.
Checking the properties is routine.\qquad
$\square$

Let $f\: X^{\updot}\to Y^{\updot}$
be a map of cosimplicial objects in
$\scrC$.
We say that $f$ is a {\it fibration} if
each $f^p$ is a fibration in $\scrC$,
and each $X^{p+1}\to
Y^{p+1}\mathop{\times}\limits_{M^p Y}
M^p X$ is a fibration in $\scrC$. Here
$M^p X\subset\mathop{\prod}\limits_{n+1}
X^n$ is the matching object of [BK, X.4.6].
Also recall from [BK, X.3.2] that
$\Tot(X^{\updot})=
\hom(\Delta^{\updot},X^{\updot})$.

\Proclaim{Lemma 1.7 {\rm [T78, 92]}}
If $f\: X^{\updot}\to Y^{\updot}$ is a
fibration of cosimplicial objects in
$\scrC$, then 
$\Tot(X^{\updot})\to\Tot(Y^{\updot})$
is a fibration in $\scrC$.
If in addition each $f^p$ is a weak
equivalence in $\scrC$ then it is also a
weak equivalence.
\finishproclaim

\proof
By [BK, X.5], axiom (SM7) for cosimplicial
objects holds.
The lemma is the case $A=\emptyset$ and
$B=\Delta^{\updot}$:
$\Tot(Y^{\updot})=
\hom(\Delta^{\updot},Y^{\updot})\times_*
*$, and $\hom(\emptyset,Y^{\updot})=*$.
\qquad $\square$

\Proclaim{Lemma 1.8 {\rm [T78, 110]}}
If $\eta\: F\Rightarrow G$ is such
that $F(K,K')\twoheadrightarrow G(K,K')$ is
a fibration in $\scrC$ for all $K$, $K'$
in $\bfK$, then $\hoend F\to\hoend G$ is
a fibration in $\scrC$.
If in addition each
$F(K,K')\twoheadSimArrow{5}
G(K,K')$ is a weak equivalence, so is
$\hoend F\to \hoend G$.

In particular, if each $F(K,K')$ is fibrant in
$\scrC$ then $\hoend F$ is fibrant in
$\scrC$.
\finishproclaim

\proof
If each $\eta_{KK'}$ is a fibration (resp.
trivial fibration) then so is each component of
$\prod\nolimits^* F\to\prod\nolimits^*G$.
Indeed $\prod^* F\to\prod^* G$ is a
fibration in $\Delta\scrC$ [BK, XI.5.3] 
[T78, 106].
Hence this follows from lemma 1.7.\qquad
$\square$

\Proclaim{Corollary 1.8.1}
The homotopy end  preserves (pointwise)
 homotopy fiber sequences.
\finishproclaim

\Proclaim{Lemma 1.9 {\rm [T82, 83] [T83, 178]}}
Suppose that $\bfK$ has an initial object
$K_0$.
Then for any functor $F\:\bfK\to\scrC$,
the natural map is a weak equivalence:
$$
\holim\limits_{\bfK} F\arrowsim
\Map(\Delta[0], F(K_0)).
$$
\finishproclaim

\proof
This is a special case of [BK, {\rm XI.4.1
(iii)}].\qquad $\square$

We conclude this section with a
generalization of the above results to a
complete ``right model category enriched
over finite simplicial sets;'' see
definitions 1.11 and 1.12 below.
The reader may easily verify that a
complete simplicial closed model category
$\scrC$ has such a structure, using lemma
1.6 and the adjunction between $\Map(-,X)$ and
$\HOM_{\scrC}(-,X)$. We will use this generality in \S7 below.

The following result is straightforward,
once one checks (see [T78, 68--92]) that (SM7)
holds, by the transfinite proof in [BK, X.5].

\Proclaim{Proposition 1.10}
If $\scrC$ is a right model category
enriched over finite simplicial sets, and
$\scrC$ is complete, then lemmas 1.6, 1.7, 1.8
and 1.9 hold in $\scrC$.
\finishproclaim

\Proclaim{Definition 1.11 {\rm [T85, 145]}} \rm
A {\it right model category} is a category
$\scrC$ with a terminal object $*$ and two
subcategories $\we(\scrC)$ and
$\fib(\scrC)$ so that: every isomorphism of
$\scrC$ is in both $\we(\scrC)$ and
$\fib(\scrC)$;  $\we(\scrC)$ is saturated
and closed under retractions; and 

\parindent=33pt
\smallskip
\Item{\sans (RM1)}
The pullback $B\times_A C \twoheadrightarrow B$ of a fibration
$C\twoheadrightarrow A$ along any map $B\to A$ exists and is a
fibration. If either $B\to A$ or $C\to A$ is in
$\we(\scrC)$, so is its pullback.
\smallskip
\Item{\sans (RM5)}
Every map $A\to B$ factors as $A\arrowsim A' \twoheadrightarrow B$,
the composite of a weak equivalence and a fibration.
\finishproclaim

Note that the subcategory $\scrC_{\fib}$ of {\it fibrant} objects $A$
(those for which $A\to*$ is a fibration) is a {\it category of fibrant
objects} in the sense of Brown [Br], with the added property that weak
equivalences are closed under retractions.

\parindent=25pt
Let $\scrS_f$ denote the category of finite simplicial sets.

\Proclaim{Definition 1.12 {\rm [T85, 146]}} \rm
We say that a complete right model
category $\scrC$ is {\it right enriched} over
$\scrS_f$ if there is a ``mapping object'' functor
$$
\Map(\,\,,\,\,)\:
\scrS_f^{\op}\times\fib(\scrC)\to\fib(\scrC)
$$
satisfying the following axioms.

\parindent=31pt
\smallskip
\Item{\sans (RE1)}
For each $K$ in $\scrS_f$:
$\Map(K,*)\cong*$; the functor $\Map(K,-)$
preserves both fibrations and base change
along fibrations in $\scrC$.

\Item{\sans (RE2)}
For each $C$ in $\scrC$:
$\Map(\emptyset,C)\cong*$; the functor
$\Map(-,C)$ sends
cofibrations in $\scrS_f$ to fibrations in
$\scrC$, and pushouts along cofibrations in
$\scrS_f$ to pullbacks along fibrations in
$\scrC$.

\Item{\sans (RE3)}
Given  a cofibration $K\maprightarrowtail{i} L$ 
in $\scrS_f$ and a fibration 
$B\maptwohead{5}{p} C$ 
in $\scrC$, 
$$
\Map(L,B)\to\Map(K,B)\times_{\Map(K,C)}
\Map(L,C)
$$
is a fibration.
It is a weak equivalence of either $i$ or
$p$ is.

\Item{\sans (RE4)}
There is a weak equivalence $\omega\:
C\arrowsim \Map(\Delta[0],C)$, natural in $C$,
so that $\partial_0\omega=\partial_1\omega\:C
\to\Map(\Delta[1],C)$.

\Item{\sans (RE5)}
For each family $\{C_i\}_{i\in I}$ of
fibrant objects in $\scrC$, and each $K$ in
$\scrS_f$, the product
$\prod C_i$ is fibrant and the canonical
morphism
$$
\Map\left(K,\prod C_i\right)\to\prod\Map
(K,C_i)
$$
is an isomorphism in $\scrC$.
\goodbreak

\parindent=25pt
\medskip
\Proclaim{Remark} \rm
Axioms (RE4) and (RE5) are not needed to prove proposition 1.10. 
They are called (RWN) and $(\aleph$EP), respectively,
in [T85, 148--9].
\finishproclaim

Here is an elementary result we shall need
for theorem 7.2.

\Proclaim{Lemma 1.13 {\rm [T85, 168]}}
Suppose given a functor $F\:\bfK\to
\fib(\scrC)$.
Then there is a functor
$(\holim F)\:\bfK\to \scrC$ defined by
$(\holim F)(K)=
\holim\limits_{(K/\bfK)}
F$, and there is a natural (pointwise)
weak equivalence
$$
F\arrowsim \holim F.
$$
\finishproclaim

\proof
For each $K$ in $\bfK$, $K/\bfK$ has an
initial object $(K$=$K)$ and a functor
$K/\bfK\to \bfK$ sending $K\to K'$ to $K'$.
By assumption, $FK\to FK'$ is a fibration
for each $K\to K'$ in $K/\bfK$.
Hence (RE4) and lemma 1.9 yield
natural weak equivalences for all $K$:
$$
FK\longmaprightsubsup{20}{\text{\rm RE4}}{\sim}
 \Map(\Delta[0], FK)
\longmaprightsubsup{20}{1.9}{\sim}
\holim\limits_{(K/\bfK)} F =(\holim\,F)(K)
F.\qquad \square
$$

\Proclaim{Remark 1.14} \rm
It is easy to see that
$$
\holim\limits_{\bfK}F=\limleft\limits_{K}
\holim\limits_{K/\bfK}F=\limleft\limits_{K}
(\holim\,F)(K).
$$
\finishproclaim

\Proclaim{Definition 1.15} \rm
We say that $\scrC$ is a {\it left model
category} if $\scrC^{\op}$ is a right
model category.
That is, $\scrC$ has an initial object
$\emptyset$ and two subcategories
$\we(\scrC)$ and $\cof(\scrC)$, both
containing all isomorphisms, so that
$\we(\scrC)$ is saturated and closed under
retractions, and 

\parindent=31pt
\smallskip\noindent
\Item{\sans (LM1)}
The pushout $P$ of a
cofibration $c\: A\rightarrowtail C$ with
any map $f\: A\to B$ exists, and $B\to P$
is a cofibration.
Moreover, if either $f$ or $c$ is a weak
equivalence, so is its pushout.

\smallskip\noindent
\Item{\sans (LM5)}
Every map $A\to B$ factors
$A\rightarrowtail B'\arrowsim B$.

\parindent=25pt
\medskip
We say that $\scrC$ is {\it left enriched} 
over $\scrS_f$ if $\scrC^{\op}$ is right
enriched by a functor
$$
\otimes\: \scrS_f\times \cof(\scrC)
\to\cof(\scrC).
$$
That is, $\otimes$ should  satisfy axioms
(LE1)--(LE5) dual to (RE1)--(RE5).
\finishproclaim

\newpage

\centerline{\sectionfont \S2 Homotopy limits of natural systems}

\bigskip\noindent
Let $\bfK$ be a category.
The {\it subdivision category} $\Sd\,\bfK$
is defined to be the category whose
objects are morphisms in $\bfK$, and whose
morphisms from $f\:A\to B$ to $f'\:A'\to B'$ are factorizations of
$f'$ through $f$:
$$\CD
B @>>> B'\\
@A{f}AA @AA{f'}A\\
A @<<< A'
\endCD$$
This construction is due to Quillen, and has many names in the literature.
Baues calls $\Sd\,\bfK$ the {\it category of factorizations} in $\bfK$
in [Baues, p. 232]. Dwyer and Kan call it the {\it Twisted arrow
category} $a\bfK$ in [DK].  (It is {\it not} the ``subdivision
category'' of [DK], or of [Mac]). 

Following Baues, a {\it natural system} $F$ on $\bfK$ with
values in $\scrC$ is a functor $F\:\Sd\,\bfK\to\scrC$.
If $\scrC$ is complete and $\bfK$ is small
we can form the products in $\scrC$:
$$
\prod\nolimits^p F=
\prod\limits_{K_0\to\cdots\to K_p}F(K_0\to K_p)\,\,.
$$

Since the indexing is over all functors $\bfp\to\bfK$,
the naturality of $F$ implies that 
$\prod\nolimits^* F$ is a cosimplicial object of $\scrC$.
Indeed, for $i\:\bfp\to \bfq$ in $\Delta$ there is a natural
map $\prod\nolimits^p F\to\prod\nolimits^q F$, 
induced by the morphism in $\Sd\,\bfK$
from $K_{i0}\to K_{ip}$ to $K_0\to K_q$.
Part of the cosimplicial identities depend upon the following observation.
Suppose that $K_0\to K_1$ and $K_p\to K_{p+1}$ are the identity.
Then for every $K_1\to K_p$:
$$
F(K_0\to K_p)=F(K_1\to K_p)=F(K_1\to
K_{p+1})\,\,.
$$

We remark that 
$\prod\nolimits^* F$ is the analogue of the chain complex used in 
[PW, 1.2] to define topological Hochshild homology; 
cf.~[JP, 3.1] and [BW].
The Baues-Wirsching cohomology [Baues, IV.5.1] of $\bfK$ with
coefficients in a natural system $F\:\Sd\,\bfK\to\bfA \bfb$ is
the cohomology of the cochain complex associated to 
$\prod\nolimits^* F$ by the Dold-Kan correspondence.
([Baues] uses the Bousfield-Kan orientation.)
\finishproclaim

\Proclaim{Definition 2.1 {\rm [T78, 114]}} \rm
Let $F\:\Sd\,\bfK\to\scrC$ be a natural
system on a small category $\bfK$, with values
in a complete category $\scrC$ equipped
with a mapping object functor.
Set
$$
\hobaw(F)= 
\Tot\left(\prod\nolimits^* F\right)\,\,,
$$
where $\Tot$ is defined in 1.1 above. The name honors Baues and
Wirsching, who introduced natural systems in [BW] as a way to study
homotopy categories.
\finishproclaim

\Proclaim{Compatibility 2.2} \rm
There is a natural functor $\Sd\,\bfK\to
\bfK^{\op}\times\bfK$ sending $A\to B$ to
$(A,B)$.
If $F$ factors as
$\Sd\,\bfK\to\bfK^{\op}\times\bfK
\longmapright{20}{F'}\scrC$ then 
$\prod\nolimits^* F$ is the same
cosimplicial object as $\prod^* F'$
in the previous
section, so $\hobaw(F)$ is the homotopy
end of $F'$:
$$
\hoend F'=\hobaw(F)\,\,.
$$
As we have seen in~1.4, the homotopy limit is
also a special case of this construction.
If $\Ftil\:\bfK^{\op}\to\scrC$ and $F$ is
the composite of $\Ftil$ with the
projection $\Sd(\bfK)\to\bfK^{\op}$ then
$$
\holim(\Ftil)=\hobaw(F)\,\,.
$$
\finishproclaim

\Proclaim{Functoriality 2.3} \rm
A natural transformation $F\Rightarrow G$
of natural systems induces a map
$\prod\nolimits^* F\to\prod\nolimits^* G$,
and so a natural map
$\hobaw(F)\to\hobaw(G)$.
The functor 
$\hobaw\:\Cat(\Sd\,\bfK,\scrC)\to\scrC$
commutes with all limits.
Given a functor $\Phi\:\bfL\to\bfK$ and
$F\:\Sd\,\bfK\to\scrC$, the diagram
$$
\CD
\Cat(\Sd\,\bfK,\scrC) 
@>{\Phi}>>
\Cat(\Sd\,\bfL,\scrC)\\
\rlap{$\raise10pt\hbox{\null\kern-10pt$\ssize\hobaw$} 
\null\kern-15pt
\vbox{\beginpicture
\setcoordinatesystem units <.50cm,.50cm>
\linethickness=1pt
\setshadesymbol ({\thinlinefont .})
\setlinear
%
%
\linethickness= 0.500pt
\setplotsymbol ({\thinlinefont .})
\plot  6.032 20.955  8.572 18.733 /
%
%
\plot  8.340 18.852  8.572 18.733  8.423 18.948 /
\linethickness=0pt
\putrectangle corners at  6.007 20.980 and  8.598 18.707
\endpicture}$}
&&\null\kern-40pt\vbox{\beginpicture
\setcoordinatesystem units <.50cm,.50cm>
\linethickness=1pt
\setshadesymbol ({\thinlinefont .})
\setlinear
%
%
\linethickness= 0.500pt
\setplotsymbol ({\thinlinefont .})
\plot 10.795 21.273  8.572 18.733 /
%
%
\plot  8.692 18.965  8.572 18.733  8.788 18.882 /
\linethickness=0pt
\putrectangle corners at  8.547 21.298 and 10.820 18.707
\endpicture}\null\kern-15pt
\raise10pt\hbox{$\ssize\hobaw$}\\
\vspace{-5pt}
&\scrC&
\endCD
$$
commutes up to a natural $2$-cell
$\varphi$ induced by the projection
$\prod\nolimits_{\bfK}^* F\to
\prod\nolimits_{\bfL}^* F$.
\finishproclaim

\goodbreak
\Proclaim{Lemma 2.4 {\rm [T78, 114] [T81, 164]}}
Let $\scrC$ be a right model category
enriched over $\scrS_f$, or a complete
closed model category.
If $F\Rightarrow G$ is such that
$F(k)\twoheadrightarrow G(k)$ is a fibration
(resp. trivial fibration) in $\scrC$ for
every $k$ in $\Sd\,\bfK$, then
$\hobaw(F)\to\hobaw(G)$ is a fibration
(resp. trivial fibration) in $\scrC$.

In particular, if each $F(k)$ is fibrant
then so is $\hobaw(F)$.
Moreover, $\hobaw$ preserves homotopy
fiber sequences.
\finishproclaim
\goodbreak

In fact, the proof of lemma 1.8 (or 1.10)
goes through in this setting.

Now consider the comma category
$k/\Sd\,\bfK$ for a fixed morphism $k$ in
$\bfK$, and write $F_k$ for the
restriction of $F\:\Sd\,\bfK\to\scrC$ to
$k/\Sd\,\bfK$.
Since $k$ is an initial object of the
comma category $k/\Sd\,\bfK$
we have a canonical weak
equivalence (as $F(k)$ is fibrant; cf.
[BK, XI.2.5]):
$$
F(k)\arrowsim
\holim\limits_{k/\Sd\,\bfK}F_k\,\,.
$$

\Proclaim{Remark {\rm [T78, 120]}} \rm
If $F\:\Sd\,\bfK\to\scrC$ takes fibrant
values, then $\hobaw(F)$ is
weak equivalent to $\holim\nolimits_{\Sd\,\bfK}
(F)$.
This result is suggested by [BW, 4.4] and [DK, 3.3].
In fact, if $F'(K\to K')$ is defined to be
$\holim\limits_{K'/\bfK}$ of 
$(K'\to L)\mapsto G(K\to
L)$, then $\hobaw (F')=\holim(F)$.
\finishproclaim

\Proclaim{Theorem 2.5 {\rm [T78, 117]}}
If $F\:\Sd\,\bfK\to\scrC$ is pointwise fibrant,
then
$$
\hobaw(F)\cong\ker\left\{
\prod\limits_{K_0}\holim\limits_{K/\Sd\,\bfK}
F \rightrightarrows
\prod\limits_{k}\holim\limits_{k/\Sd\,\bfK} F_k
\right\}\,\,.
$$
Since $\hobaw$ preserves homotopy
equivalences of such fibrant valued $F$,
this shows that $\hobaw(F)$ is the right
derived functor of
$$
F\longmapsto \ker\left\{
\prod\limits_{K_0}F(K_0=K_0)\rightrightarrows
\prod\limits_{K_0\to K_1}F(K_0\to
K_1)\right\}
$$
for $F\:\Sd\,\bfK\to\scrC$.
\finishproclaim

\proof
The right side is
$$
\align
\ker\prod\nolimits_{K}\Biggl\{\Tot
\left(
\mathop{\prod\nolimits^*}\limits_{K/\Sd\,\bfK}
F\right) &\rightrightarrows
\prod\nolimits_{k}\Tot\left(
\mathop{\prod\nolimits^*}\limits_{k/\Sd\,\bfK}
F_k\right)\\
\vspace{15pt}
=\Tot\Biggl\{\ker\prod\nolimits_{K}
\left(\mathop{\prod\nolimits^*}
\limits_{K/\Sd\,\bfK}
F\right) &\rightrightarrows
\prod\limits_{k}\left(
\mathop{\prod\nolimits^*}\limits_{k/\Sd\,\bfK}
F_k\right)\Biggr\}\,\,.
\endalign
$$
This is $\Tot$ of a cosimplicial object in
$\scrC$ which in degree $p$ is an
equalizer of 
$\prod\nolimits_{K}\prod\nolimits_{\lambda}
F(\lambda_p)\rightrightarrows 
\prod\nolimits_{k}\prod\nolimits_{\mu}
F(\mu_p)$, where $\lambda$ and $\mu$ run
over diagrams in $K/\Sd\,\bfK$ and
$k/\Sd\,\bfK$, respectively.
This equalizer is easily seen to be $
\prod\nolimits^p F$, so the cosimplicial
object is just the simplicial
replacement $\prod\nolimits^* F$.
We are done, since $\hobaw(F)$ is defined
to be $\Tot\left(\prod^* F\right)$.\qquad
$\square$

Thomason defines an even more exotic generalization
of a homotopy limit, which we christen $\hofoo$.
To define it, recall that the objects of the comma
category $\Delta/\bfK$ are pairs
$(p,K\:\bfp\to\bfK)$, i.e., diagrams
$K_0\to K_1\to\cdots\to K_p$ in $\bfK$,
while a morphism $(p,K)\to(q,L)$ consists
of a simplicial map $\sigma\: \bfp\to\bfq$
such that $K=L\sigma$.
It will be convenient to abbreviate
$(p,K)$ as $K$ and write $\#K$ for the
first entry  $p$ of $(p,K)$.

Given a functor $F\: \Delta/\bfK\to\scrC$,
let $\prod\nolimits^{*}F$ be the
cosimplicial replacement defined by
$$
\prod\nolimits^p
F=\prod\limits_{\#K=p}F(K)\,\,.
$$
If $\sigma\:\bfq\to\bfp$ is a morphism in
$\Delta$, we define
$\sigma_*\:\prod\nolimits^q
F\to\prod\nolimits^p F$ by requiring that
the component corresponding to
$K\:\bfp\to\bfK$ is the projection onto
the component $F(K\sigma)$ followed by the
map $F(K\sigma)\to F(K)$ associated to the
morphism $K\sigma\to K$ induced by $\sigma$.
$$
\bfq\longmapright{20}{\sigma}\bfp
\longmapright{20}{K}\bfK
$$
\goodbreak
\Proclaim{Definition 2.6 {\rm [T85, 36]}} \rm
Let $F\:\Delta/\bfK\to\scrC$ be a functor,
where $\scrC$ is a complete category
equipped with a mapping object functor.
Then
$$
\hofoo(F)=
\Tot\left(\prod\nolimits^{*}F\right)
=\int\limits_{p\in\Delta}
\Map(\Delta[p],\prod\limits_{\#K=p}F(K))
\,\,.
$$
\finishproclaim
\goodbreak

\Proclaim{Lemma 2.7 {\rm [T83, 158]}}
There is a functor
$(K_1,K_2)\mapsto\Map(\Delta[\#K_1],FK_2)$
from
$(\Delta/\bfK)^{\op}\times
(\Delta/\bfK)$ to $\scrC$, and $\hofoo(F)=
\hofoo_{\Delta/\bfK}(F)$ is its end:
$$
\hofoo(F)\cong\int\limits_{K\in\Delta/\bfK}
\Map(\Delta[\#K],FK)\,\,.
$$
\finishproclaim

\proof
It suffices to check that the two ends
have the same universal mapping property
with respect to objects of $\scrC$.
The details are routine, and duly verified in
[T83, 159].\qquad $\square$

There are forgetful functors
$$
\gather
(\Delta/\bfK)\longrightarrow \Sd\,\bfK
\longrightarrow\bfK^{\op}\times 
\bfK\longrightarrow \bfK\\
(K_0\to\cdots\to K_p)\longmapsto(K_0\to
K_p) \longmapsto (K_0,K_p)\longmapsto
K_0\,\,.
\endgather
$$
\goodbreak

\Proclaim{Proposition 2.8 {\rm [T85, 152]}}
For each $F\:\Delta/\bfK\to\scrC$:
\smallskip
\Item{\sans (a)}
If $F$ factors through
$F'\:\Sd\,\bfK \to\scrC$ then $\hofoo(F)=\hobaw(F')$;
\smallskip
\Item{\sans (b)}
If $F$ factors through $F''\:\bfK^{\op}\times\bfK\to\scrC$ then
$\hofoo(F)=\hoend F''$;
\smallskip
\Item{\sans (c)}
If $F$ factors through $F'''\:\bfK^{\op}\to\scrC$ then
$\hofoo(F)=\holim(F''')$.
\finishproclaim

Indeed, (a) is routine since $\prod^* F$ is the same cosimplicial
object as $\prod^* F'$. Parts (b) and (c) follow from Compatibility~2.2.

\Proclaim{Remark 2.8.1} \rm
$\Delta/(\bfK^{\op})$ is not
$(\Delta/\bfK)^{\op}$, and
$\Sd(\bfK^{\op})$ is not $(\Sd\,\bfK)^{\op}$.
\finishproclaim

Here is the generalization of lemma 1.8;
again the proof of 1.8 goes through.

\Proclaim{Theorem 2.9 {\rm [T85, 155]}}
Let $\scrC$ be a complete right model
category enriched over $\scrS_f$.
Suppose in addition that towers and
products preserve fibrations and
trivial fibrations in $\scrC$.
Then for every small category $\bfK$ and
every pointwise fibrant functor
$$
F\:\Delta/\bfK\to\scrC\,\,,
$$
the object
$\hofoo(F)$ is fibrant in $\scrC$.
If $\eta\:F\to G$ is a pointwise fibration
between pointwise fibration functors then
$\hofoo(\eta)\:\hofoo(F)\to\hofoo(G)$ is a
fibration.
If in addition $\eta$ is a weak
equivalence, so is $\hofoo(\eta)$.
\finishproclaim

\Proclaim{Variant 2.10 {\rm [T85, 6]}}\rm
Given $F\:(\Delta/\bfK)^{\op}\to
\scrC_{\cof}$, we have a simplicial object
in $\scrC_{\cof}$ whose $p$th
term is $\coprod\limits_{\#K=p}F(K)$.
If $\scrC$ is cocomplete and
left enriched (1.15), one can
define $\hocofoo(F)$ as the coend
$\int^{\Delta}\Delta [\#K]\otimes
\coprod\nolimits_{*}F$.
The dual assertions for 2.7--2.8 are
explored in [T85].
The dual of theorem 2.9 is given in [T81, 10, 159]
and on pp.~46--48 of [T86].
\finishproclaim

\newpage

\centerline{\sectionfont \S3 Thomason's
Axioms for homotopy theory}

\bigskip\noindent
Let $\scrC$ be a category.
If $W$ is a distinguished family of
morphisms, closed under composition and
containing all identity maps, 
it is convenient to think of $W$ as the
morphisms of subcategory $\scrW$ having
the same objects as $\scrC$.
We say that $\scrW$ is {\it saturated} if
given any pair $(f,g)$ of composable maps
in $\scrC$, whenever two of $f$, $g$ and
$fg$ are in $\scrW$ so is the third.

Let $\scrC$ be a category equipped with
three distinguished subcategories
$\cof(\scrC)$,  $\fib(\scrC)$ and
$\we(\scrC)$ having the same objects as
$\scrC$.
Morphisms in these subcategories will be
called {\it cofibrations} (written
$\rightarrowtail$), {\it fibrations}
(written $\twoheadrightarrow$) and {\it
weak equivalences} (written $\arrowsim$),
respectively.

A {\it trivial cofibration}, or 
{\it equicofibration}, (written 
$\longmaprightarrowtail{\sim}$) is a map which
is both a cofibration and a weak
equivalence; a {\it trivial fibration}, or
{\it equifibration},
(written $\maptwohead{10}{\sim}$) is a map
which is both a fibration and a weak
equivalence.

We suppose that $\scrC$ has an initial
object $\emptyset$ and a terminal object
$*$.
Consider the following self-dual axioms
[T78, 156].

\parindent=35pt
\Item{(TM0)}
Every isomorphism of $\scrC$ is in
$\cof(\scrC)$, $\fib(\scrC)$ and
$\we(\scrC)$.

\Item{(TM1)}
The pushout $P$ of a cofibration $c\:
A\rightarrowtail C$ with any map $f\: A\to
B$ exists, and $B\to P$ is a cofibration.
Moreover, if either $f$ or $c$ is a weak
equivalence, so is its pushout.
The dual assertion (RM1) for the pullback
$B\times_A C$ of a fibration 
$C\twoheadrightarrow A$ and any map $B\to
A$ must also hold.

\Item{(TM2)}
(Quillen's Axiom CM2)
The subcategory $\we(\scrC)$ is saturated.

\Item{(TM3)}
(Quillen's Axiom CM3w)\footnote"$^\dagger$"{Axiom
 (TM3) is not needed for cofibrations or
fibrations, only for weak equivalences.
See [T79,~26] and [T83,~152].}
Any retract of a weak equivalence is a
weak equivalence.

\Item{(TM4)}
(Quillen's Lifting Axiom CM4)
Suppose given a commutative square
$$
\CD
A @>>> X\\
\Ldownarrowtail{i} @. \Rtwoheadmapdown{p}\\
B @>>> Y
\endCD
$$
in which $i$ is a cofibration and $p$ is a
fibration.
Then a map $B\to X$ exists, factoring both
$A\to X$ and $B\to Y$, provided that
either $i$ or $p$ is a weak equivalence.
Any such map $B\to X$ will be called a 
{\it fill-in}, or {\it factorization}.

\medskip
\Item{(CM5)}
Any map $A\mapright{f} B$ factors as $A
\rightarrowtail B'\arrowsim B$, and also
as $A\arrowsim A'\twoheadrightarrow B$.

\Item{(TM5)}
The following factorizations of $A\arrowsim B$
can be made functorial in $f$:
$$
\align
&A\rightarrowtail B'\arrowsim
  B\tag\hbox{c}\\
&A\arrowsim A'\twoheadrightarrow
B\tag\hbox{f}
\endalign
$$

\medskip
\parindent=25pt
The functoriality in (TM5) may be
expressed as follows.
Let $\scrC^{\bfone}$ denote the category
of arrows in $\scrC$, i.e., functors
$\bfone\mapright{f}\scrC$, and
$\scrC^{\bftwo}$ the category of composable pairs
of arrows in $\scrC$, i.e., functors
$\bftwo\to\scrC$.
Then there are functors
$$
T,M\: \scrC^{\bfone}\to\scrC^{\bftwo}
$$
sending the arrow $f\:A\to B$ in
$\scrC^{\bfone}$ to the
composable pairs (with composition $f$)
$$
A\rightarrowtail T(f)\arrowsim B\qquad
\text{and}\qquad A\arrowsim
M(f)\twoheadrightarrow B
$$
respectively. 
See [T79, 72-73] and [T78, 161].

We call an object $C$ {\it cofibrant} if
the canonical map $\emptyset\to C$ is in
$\cof(\scrC)$; we call $C$ {\it fibrant}
if the canonical map $C\to *$ is in
$\fib(\scrC)$.

\Proclaim{Definition 3.1} \rm
A {\it basic model category} is a
category $\scrC$, equipped with
distinguished subcategories $\cof(\scrC)$,
$\fib(\scrC)$ and $\we(\scrC)$, an initial
object $\emptyset$ and a terminal object $*$,
satisfying axioms (TM0)-(TM4) and (CM5).
If in addition it satisfies axiom (TM5),
we call $\scrC$ a {\it Thomason model
category}.

Comparing definitions, we see that any
basic model category is both a left model
category and a right model category.
\finishproclaim

The following consequence of the axioms is
of fundamental importance, so we include
it here.

\Proclaim{Glueing Lemma 3.2}
Suppose that $\scrC$ satisfies axiom
(TM0), (TM1) and (TM2).
Given a commutative diagram
$$
\CD
B @<<< A &\longrightarrowtail &C\\
@V{\sim}VV @VV{\sim}V @VV{\sim}V\\
B' @<<< A' &\longrightarrowtail &C'
\endCD
$$
in which the vertical maps are weak
equivalences, the pushout $B\sqcup_A
C\to B'\sqcup_{A'} C'$ is a weak
equivalence.
\finishproclaim

\proof
By (TM1) the map $C\to A'\sqcup_A C$ exists
and is a weak equivalence.
Its composition with $g\: A'\sqcup_A C\to$ 
is $C\arrowsim C'$, so $g$ is a weak
equivalence by (TM2).
But then we may use (TM0) and (TM1) twice
to get
$$
B\sqcup_A C\arrowsim B'\sqcup_A
C\cong B'\sqcup_{A'}A'
\sqcup_A C\arrowsim B'\sqcup_{A'}C'\,\,.
$$
The composition is a weak equivalence, as
desired.\qquad $\square$

\Proclaim{Remark 3.2.1} \rm
(TM0) is only used for $\we(\scrC)$ and
(TM1) is only used for cofibrations.
\finishproclaim

\newpage

\centerline{\sectionfont \S4 Examples}

\medskip
As usual, the paradigm of a Thomason model
category is the category of
topological spaces and continuous maps,
where the fibrations are Serre fibrations
and weak equivalences are homotopy
equivalences; cofibrations are defined by
the lifting property (TM4).

\Proclaim{Definition 4.1} \rm
More generally, any proper model category
in the sense of Quillen [QH] [QR] [BF]
[GJ] is a basic model category; the
adjective {\it proper} describes the hard part of axiom (TM1).
The easy part of (TM1) is axiom (M4) in
Quillen's original definition of a model
category [QH]; see [GJ, II.1.3].
\finishproclaim

I do not know of any cocomplete proper model category
in which the factorizations required by
axiom (CM5) are not functorial.
We will see in 4.8 below that if $\scrC$ admits a
{\it small object argument} then $\scrC$
satisfies (TM5) and is a Thomason model category.
Compare [GJ, I.9.2].

Here is another difference between
Quillen's and Thomason's axioms.
Thomason's Factorization Axiom (TM5)
allows us to factor $f$ as $hi$, where $i$
is a cofibration and $h$ is a weak
equivalence; in Quillen's Factorization
Axiom (CM5) $h$ must also be a fibration.
Therefore if $f$ is a map having the left
lifting property with respect to all
fibrations which are equimorphisms,
Quillen's axioms (CM3) and (CM5) imply that
$f$ must be a cofibration.
This need not be the case in a Thomason
model category, as the following example shows.

\Proclaim{Example 4.2} \rm
Here is an example of a Thomason model
category which is not a Quillen model
category.
Let $\Ch =\Ch^b(R)$ denote the category of
bounded chain complexes of modules over a
ring $R$, with weak equivalences being
quasi-isomorph-\break
isms (maps inducing isomorphisms on homology).

Let $\cof(\Ch)$ denote the category of
injections $A\hookrightarrow B$ 
whose cokernel $B/A$ is a chain
complex of projective modules, and let
$\fib(\Ch)$ denote the category of
surjections whose kernel is a chain
complex of injective modules.
If $R$ is a noetherian regular ring of
finite global dimension, then $\Ch^b(R)$
is a Thomason model category but not a
Quillen model category 
(unless $R$ is semisimple).
\goodbreak

Of course, $\Ch^b(R)$ has the same homotopy
theory as the Quillen model structure of
[QH, I.1.2], in which every surjection is
a fibration.
See also [QH, I.4.12].
\finishproclaim

\Proclaim{Example 4.3} \rm Baues' {\it cofibration category} is
another context in which one can do
homotopy theory; this theory is developed in [Baues].
By definition, a cofibration category is a
category $\scrC$ equipped with two classes
of morphisms ($\cof$ and $\we$) such that
Baues' axioms (C1)--(C4) hold.
Thus Baues' cofibration categories lack
all fibrations present in a Thomason model category.
\finishproclaim

Baues' axioms (C1), (C2), and (C3) for
$(\scrC,\cof,\we)$ are axioms (TM0) ,
(TM1), (TM2) and (CM5c) for cofibrations and
weak equivalences.
Baues calls an object $X$ a {\it fibrant
model} if each trivial cofibration 
$X\maprightarrowtail{\sim} Y$ 
admits a retraction, and Baues' final axiom is:

\smallskip
\Item{\rm (C4)}
For each object $A$ in $\scrC$ there is a
trivial cofibration $A\maprightarrowtail{\sim} X$
where $X$ is a fibrant model.

\medskip\noindent
If $\scrC$ is Thomason model category,
then $(\scrC,\cof,\we)$ is a left model
category in the sense of 1.13 above.
However, $(\scrC,\cof,\we)$ may not satisfy
Baues' axiom (C4), so $\scrC$ may not be a
cofibration category in the sense of
[Baues].

For purposes of comparison, recall 
from [Baues, I.2.6]
that if $(\scrC,\cof,\fib,\we)$ is a model
category in the sense of Quillen [QH],
then the subcategory $\scrC_{\cof}$ of all
cofibrant objects in $\scrC$ is a
cofibration category.

\Proclaim{Example 4.4} \rm
Recall from [TT, 1.2.4] that a {\it
biWaldhausen category} is a category
$\scrC$ with a zero object $0$, together
with subcategories $\cof(\scrC)$,
$\fib(\scrC)$ and $\we(\scrC)$ so that:

\medskip\noindent
\Item{\rm (i)}
$(\scrC,\cof,\we)$ and
$(\scrC^{\op},\fib,\we)$ are Waldhausen
categories [TT, 1.2.3];

\Item{\rm (ii)}
the canonical map $A\sqcup B\to A\times B$
is an isomorphism; and

\Item{\rm (iii)}
the cofibration sequences in $\scrC$ and
$\scrC^{\op}$ are dual to each other.
\finishproclaim

\medskip\noindent
Note that every object  of $\scrC$ is both
cofibrant and fibrant by axiom 1.2.1.2 of
[TT].
A biWaldhausen category satisfies (TM0)
and (TM1) by axioms 1.2.1.1, 1.2.1.3,
1.2.3.1 and 1.2.3.2 of [TT].
We say that $\scrC$ is {\it saturated} if
axiom (TM2) holds, i.e., if $\we(\scrC)$
is a saturated subcategory.

If $\scrC$ has mapping cylinders
satisfying the cylinder axiom, and mapping
path spaces satisfying the path space
axiom (see [TT, 1.3.1 and 1.3.2]) then
$\scrC$ satisfies (TM4) and (TM5).
Thus only the retraction
axiom (TM3) is needed to make
$\scrC$ a Thomason model category.
Here is a partial converse.

\Proclaim{Lemma 4.5}
If $\scrC$ satisfies (TM0) -- (TM3),
$A\sqcup B\simeq A\times B$ for all $A$ and
$B$, and every object is both fibrant and
cofibrant, then $\scrC$ is a saturated
biWaldhausen category.
\finishproclaim

\proof
This is straightforward; since objects are
fibrant and cofibrant, axiom (i) follows from
(TM0), (TM1) and the gluing lemma 3.2.\qquad
$\square$

We caution the reader that the
functoriality axiom (TM5) need not give a
Waldhausen cylinder functor $T$, 
because the map
$A\sqcup B\rightarrowtail T(f)$ need not
be in $\cof(\scrC)$.

We now give a set-theoretic condition
under which a proper closed Quillen model category
is a Thomason model category.
It is based on II.3.4 of [QH].

\Proclaim{Definition 4.6} \rm
Suppose that $\scrC$ is a closed model
category.
We say that a set of trivial cofibrations
$\left\{A_\alpha\maprightarrowtail{
\sim} B_\alpha\right\}$ 
is a set of {\it test
cofibrations} if for each map $E\mapright{p} F$, 
$p$ is a fibration if and only if every
diagram of the following form admits a
factorization $B_\alpha\to E$.
$$
\CD
A_\alpha @>>> E\\
\downarrowtail @. @VV{p}V\\
B_\alpha @>>> F
\endCD
$$

For example, the {\it anodyne extensions}
form a set of test cofibrations for the
category $\scrC$ of simplicial sets; see
[GJ, I.4.3].
The horns 
$\left\{\Lambda_k^n\maprightarrowtail{\sim}
\Delta^n\right\}$ form a countable set of
test fibrations; see [GJ, p.~11].
\finishproclaim

Let $\kappa$ be an infinite cardinal
number.
A poset is called {\it $\kappa$-filtering}
if any subset of cardinality $<\kappa$ has
an upper bound.

\Proclaim{Definition 4.7} \rm
An object $A$ of a category $\scrC$ is
called {\it $\kappa$-small} if for every
$\kappa$-filtering directed poset $\{C_i\}$
for which $\limright\,C_i$ exists we have
$$
\limright\,\Hom_{\scrC}(A,C_i)\cong
\Hom_{\scrC}(A,\limright\,C_i)\,\,.
$$
\finishproclaim

\Proclaim{Theorem 4.8 {\rm (Small Object
Argument)}}
Let $\scrC$ be a closed model category
with arbitrary coproducts.
Suppose $\scrC$ has a set
$\{A_\alpha\maprightarrowtail{\sim}
B_\alpha\}$ of test cofibrations in which
each $A_\alpha$ is $\kappa$-small for some
fixed $\kappa$.
Then $\scrC$ satisfies (TM5).
If in addition $\scrC$ is a proper closed
model category, then $\scrC$ is a Thomason
model category.
\finishproclaim

\proof
We construct $F_\mu\: 
\scrC^{\bfone}\to\scrC^{\bftwo}$ by
induction on $\mu\le\kappa$.
Fix $f\: X\to Y$.
The functor $F_0$ sends $f$ to
$X\maprightarrowtail{\id} X\longrightarrow Y$.
We define $T_{\mu+1}$ by
$$
T_{\mu+1}(f)=\left(\coprod\limits_{S}
B_\alpha\right)\coprod\limits_{\left(\coprod
\nolimits_{S}A_\alpha\right)}T_\mu(f)
$$
where $S$ runs over all test diagrams
$$
\CD
A_\alpha @>>> T_\mu(f)\\
\downarrowtail @. \twoheaddownarrow\\
B_\alpha @>>> Y\,.
\endCD
\tag$*$
$$
The map $T_\mu\to T_{\mu+1}$ is a trivial
cofibration as in Lemma 3 of [QH, II.3.3]
and the maps $B_\alpha \to Y$ and
$T_\mu\twoheadrightarrow Y$ induce a
fibration $T_{\mu+1}\to Y$.
We define $F_{\mu+1}$ by
$$
X \maprightarrowtail{\sim} T_\mu
\maprightarrowtail{\sim} T_{\mu+1}
\twoheadrightarrow Y\,\,.
$$
For a limit ordinal $\mu$, we observe that
$X\rightarrowtail \limright\,T_\nu=T_\mu$ is a
trivial cofibration.
The map $T_\mu\to Y$ is a fibration
because for each test diagram $(*)$ there
is a $\nu<\mu$ so that $A\to T_\mu$
factors through a map $A\to T_\nu$,
$\nu<\mu$ ($A_\alpha$ is $\mu$-small as
$\mu<\kappa$), and
the factorization $B\to T_{\nu+1}\to
T_\mu$ of $(*)$ exists by the construction
of $T_{\nu+1}$.\qquad
$\square$

\newpage

\centerline{\sectionfont \S5 The Homotopy
Category}

\bigskip
Suppose that $\scrC$ is a basic model category.
The basic idea needed to construct the
homotopy category $\Ho\scrC$ is quite
simple, and goes back to Quillen [QH].
Given $A$ cofibrant and $B$ fibrant, we show in 5.4 that there
is a good equivalence relation $\simeq$ on
$\Hom_{\scrC}(A,B)$ and $[A,B]$ is the set
of equivalence classes.
Given general $A$ and $B$, axiom (CM5)
gives factorizations $\emptyset
\rightarrowtail A'\arrowsim A$ and
$B\arrowsim B'\twoheadrightarrow *$;
we set $[A,B]=[A',B']$.
As usual, there are lots of technical
details to check.

By a {\it cylinder object} for an object
$A$ of $\scrC$, we mean an object $C$ of
$\scrC$, together with a factorization
$A\sqcup A\rightarrowtail C\arrowsim A$ of
the fold map $\nabla\: A\sqcup A\to A$ (we
assume $A\sqcup A$ exists).
Dually, a {\it path object} for $B$ is an
object $M$ together with a factorization
$B\arrowsim M\twoheadrightarrow B\times B$
of the diagonal $\Delta$ (we assume
$B\times B$ exists).

\Proclaim{Definition 5.1} \rm
For $f,g\in \Hom_{\scrC}(A,B)$, define $f\simeq_L g$ to hold if the
map $f\sqcup g\:A\sqcup A\to B$ factors as 
$A\sqcup A \maprightarrowtail{\partial} C\mapright{H} B$,
where there is a weak equivalence
$q\: C\arrowsim A$ such that $q\partial$
is the fold map, i.e., $C$ is a cylinder
object for $A$.
\finishproclaim

\Proclaim{Lemma 5.2 {\rm [T78, 166]}}
If $A$ is cofibrant then $\simeq_L$ is an
equivalence relation on the set
$\Hom(A,B)$.
We define $\pi(A,B)$ to be
$\Hom(A,B)_{\simeq_L}$.
\finishproclaim

If $A$ is cofibrant and $B\to B'$
arbitrary, it is easy to see that
$\Hom(A,B)\to\Hom(A,B')$ induces a
quotient map $\pi(A,B)\to\pi(A,B')$.

\proof
(Quillen [QH, I.1.8])
Axiom (CM5) implies that the fold map
$\nabla$ factors as $A\sqcup
A\rightarrowtail C\maptwohead{10}{\sim}
A$; the composition with $f\: A\to B$ is
clearly $f\sqcup f$, so $f\simeq_L f$.
To see that the relation is symmetric,
precompose the factorization of $f\sqcup
g$ with the permutation $\tau$; the
composition
$$
A\sqcup A\mapright{\tau} A\sqcup 
A\longmaprightarrowtail{\partial}
C\mapright{H}
B
$$
will then be $g\sqcup f$, and
$q\partial\tau=\nabla\tau=\nabla$.
Finally, to see that $\simeq_L$ is
transitive, suppose given factorizations
of $f\sqcup g$ and $g\sqcup h$:
$$
A\sqcup A\rightarrowtail
C_1\longrightarrow B\qquad\text{and}\qquad
A\sqcup A\rightarrowtail
C_2\longrightarrow B\,\,.
$$
Since $A$ is cofibrant, $A=\emptyset\sqcup
A\rightarrowtail A\sqcup A\rightarrowtail
C_1$ and $A=A\sqcup\emptyset
\rightarrowtail A\sqcup A\rightarrowtail
C_2$ are cofibrations.
Hence we may apply the Glueing Lemma to
$$
\minCDarrowwidth{15pt}
\matrix
C_1 &\longleftarrowtail &A 
&\longrightarrowtail &C_2\\
\lmapdown{q_1} &&\null\,\,\rmapdown{=}
&&\null\,\,\rmapdown{q_2}\\
A @= A @= A\,.
\endmatrix
$$
This shows that the pushout $a\: C_1\sqcup_A
C_2\to A$ is a weak equivalence.
Moreover, the composites $A\rightarrowtail
C_i\to B$ are both $g$, so they induce a
map $h_1\sqcup_A h_2\: C_1\sqcup_A C_2\to B$.
If we let $\partial$ be the map
$$
A\sqcup
A=(A\sqcup\emptyset)\sqcup(\emptyset\sqcup A)
\rightarrowtail C_1\sqcup C_2 \to C_1
\sqcup_A C_2
$$
then clearly $q\partial=\nabla$ and
$(h_1\sqcup_A h_2)\circ\partial=f\sqcup h$.
Finally, $\partial$ is a cofibration
because it is the composition of two
cofibrations:
$\partial_0\sqcup A\: A\sqcup
A\rightarrowtail C_1\sqcup A$, the pushout
of $A\rightarrowtail C_1$ by $A=A\sqcup
\emptyset\rightarrowtail A\sqcup A$, and
the map $C_1\sqcup A\rightarrowtail
C_1\sqcup_A C_2$, which is the pushout of $A\sqcup
A\rightarrowtail C_2$ by $\partial_0\sqcup
A\: A\sqcup A\rightarrowtail 
C_1\sqcup A$.\qquad $\square$
$$
\CD
A @>>> A\sqcup A 
  &\longrightarrowtail &C_2\\
\downarrowtail @. 
  \Rdownarrowtail{\partial_0\sqcup A}
@.\downarrowtail\\
C_1 @>>> C_1\sqcup A 
  &\longrightarrowtail
&C_1\sqcup_A C_2
\endCD
$$
If $A$ is cofibrant, the following
argument shows that we can fix the choice
of cylinder object in the definition of
$\simeq_L$.
This fact allows us to simplify some
set-theoretic considerations later on.

\Proclaim{Lemma 5.3}
If $A$ is cofibrant, $B$ is fibrant and
$C$ is a cylinder object for $A$, then
whenever two maps $f,g\:A\to B$ satisfy
$f\simeq_L g$ there is a map $C\to B$ so
that $A\sqcup A\rightarrowtail C\to B$ is
$f\sqcup g$.
\finishproclaim

\proof [T83, 77]
We are given a factorization $A\sqcup
A\rightarrowtail C' \mapright{H'} B$ of
$f\sqcup g$, where $C'$ is a different
cylinder object for $A$.
Let $C''$ denote the pushout of $C$ and
$C'$ along $A\sqcup A$, and factor the
pushout $C''\to A$ as $C''\rightarrowtail
C'''\arrowsim A$.
Axiom (TM2) shows that the composite
$C'\rightarrowtail C''\rightarrowtail
C'''$ is a trivial cofibration.
Since $B$ is fibrant, axiom (TM4) applied
to the right side of
$$
\matrix
A\sqcup A &\longrightarrowtail
  &C' &\longmapright{20}{H'} &B\\
\downarrowtail &&\Ldownarrowtail{\sim} 
  &\lower5pt\hbox{\vbox{\beginpicture
\setcoordinatesystem units <.40cm,.40cm>
\linethickness=1pt
\setlinear
%
%
\linethickness= 0.500pt
\setplotsymbol ({\thinlinefont .})
\setdashes < 0.1270cm>
\plot  8.604 16.796 10.890 18.161 /
\setsolid
%
%
\plot 10.705 17.976 10.890 18.161 10.640 18.085 /
\setdashes < 0.1270cm>
\linethickness=0pt
\putrectangle corners at  8.579 18.186 and 10.916 16.770
\endpicture}} &\twoheaddownarrow\\
C &\longrightarrowtail &C'''
&\longmapright{20}{} &*
\endmatrix
$$
shows that $H'$ factors through $C'''$.
But then $f\sqcup g$ factors as
$$
A\sqcup A\rightarrowtail C\longrightarrow
C''' \longrightarrow B\,\,. \qquad
\square
$$

Dually, if $B$ is fibrant, we get an
equivalence relation $\simeq_R$ on
$\Hom_{\scrC}(A,B)$ by declaring
$f\simeq_R g$ if there is a diagram:
$$
\matrix
A &\longmapright{20}{} &M &\longmapleft{20}{\sim}
  &B\\
\noalign{\medskip}
\raise6pt\hbox{$_{f\times g}$}
\null\kern-5pt\rlap{$\vbox{\beginpicture
\setcoordinatesystem units <.35cm,.35cm>
\linethickness=1pt
\setshadesymbol ({\thinlinefont .})
\setlinear
%
%
\linethickness= 0.500pt
\setplotsymbol ({\thinlinefont .})
\plot  6.032 20.955  8.572 18.733 /
%
%
\plot  8.340 18.852  8.572 18.733  8.423 18.948 /
\linethickness=0pt
\putrectangle corners at  6.007 20.980 and  8.598 18.707
\endpicture}$} &&\twoheaddownarrow 
&&\llap{$\vbox{\beginpicture
\setcoordinatesystem units <.35cm,.35cm>
\linethickness=1pt
\setshadesymbol ({\thinlinefont .})
\setlinear
%
%
\linethickness= 0.500pt
\setplotsymbol ({\thinlinefont .})
\plot 10.795 21.273  8.572 18.733 /
%
%
\plot  8.692 18.965  8.572 18.733  8.788 18.882 /
\linethickness=0pt
\putrectangle corners at  8.547 21.298 and 10.820 18.707
\endpicture}$}\null\kern-5pt
  \raise6pt\hbox{$\ssize \Delta$}\\
\noalign{\medskip}
&&B\times B &
\endmatrix
$$
\Proclaim{Lemma 5.4}
If $A$ is cofibrant and $B$ fibrant, then
for all $f,g\in\Hom_{\scrC}(A,B)$ we have:
$f\simeq_L g$ if and only if $f\simeq_R
g$.
\finishproclaim

\proof
Suppose that $f\simeq_L g$ via $A\sqcup
A\rightarrowtail C\mapright{h} B$.
By (CM5) we can factor $\Delta$ as
$B\arrowsim M\twoheadrightarrow B\times
B$.
By (TM4), there is a fill-in $\varphi\:C\to M$
for the diagram:
$$
\matrix
A &\mapright{f}\,\, B\,\,\mapright{\sim} &M\\
\Ldownarrowtail{\partial_0} &&\twoheaddownarrow\\
C &\longmaprightsub{40}{(fg,h)} &B\times B
\endmatrix
$$
and $\varphi\partial_1\:A\to M$ is the map
needed to show that $f\simeq_R g$.
By duality, $f\simeq_R g$ implies
$f\simeq_L g$.\qquad $\square$

\Proclaim{Proposition 5.4 {\rm [T78,~170]}}
Let $\scrC$ be a basic model category.
If $A$ and $A'$ are cofibrant and $B$ is
fibrant, then for every weak equivalence
$w\:A'\to A$ there is a bijection
$$
w^*\:\pi(A,B)\to\pi(A',B)
$$
\finishproclaim

\Proclaim{Remark} \rm
This is [QH, I.1.10] if $w$ is a trivial
cofibration.
\finishproclaim

\proof
The map $w^*$ is well defined since by
definition $f\simeq_R g$ implies
$fw\simeq_R gw$.
To see that $w^*$ is injective, suppose given
$f,g\:A\to B$ so that $fw\simeq_L gw$.
That is, $fw\sqcup gw$ factors as
$A'\sqcup A'\rightarrowtail C\to B$.
Form the pushout $P$ of $C$ and $A\sqcup
A$ along $A'\sqcup A'$; the map $f\sqcup
g$ and $C\to B$ agree on $A'\sqcup A'$
and induce a map $G\:P\to B$.
$$
\matrix
A'\sqcup A' &\longmaprightarrowtail{\partial'}
 &C &\longmaprightsubsup{20}{q'}{\sim} &A'\\
\lrmapdown{w\sqcup w}{\sim}
&&\Big\downarrow
  &&\rmapdown{w}\\
\vspace{5pt}
A\sqcup A &\longmaprightarrowtail{\partial}
 &P &\mapright{q} &A\\
\vspace{5pt}
\null\qquad\qquad
\hbox{$\ssize f\sqcup g$}\kern-20pt
\rlap{$\lower8pt\hbox{$\vbox{\beginpicture
\setcoordinatesystem units <.60cm,.60cm>
\linethickness=1pt
\setshadesymbol ({\thinlinefont .})
\setlinear
%
%
\linethickness= 0.500pt
\setplotsymbol ({\thinlinefont .})
\plot  6.032 20.955  8.572 18.733 /
%
%
\plot  8.340 18.852  8.572 18.733  8.423 18.948 /
\linethickness=0pt
\putrectangle corners at  6.007 20.980 and  8.598 18.707
\endpicture}$}$} 
 &&\raise10pt\hbox{$\null\,\,\,\rmapdown{G}$} &&\\
\vspace{5pt}
&&B &&
\endmatrix
$$
Similarly $\nabla$ and $wq'$ induce a map
$q\:P\to A$ so that $q\partial=\nabla$.
Since $A$ and $A'$ are cofibrant, the
glueing lemma shows that $w\sqcup w\: A'
\sqcup A'\to A\sqcup A$ is a weak
equivalence.
By axiom (TM1), its pushout $C\arrowsim P$
is a weak equivalence.
Its composition with $q$ is a weak
equivalence, so $q$ is a weak equivalence
by (TM2).
Thus $P$ gives $f\simeq_L g$.

To see that $w^*$ is onto we have to lift
any $A'\mapright{f} B$ to $A$.
If $w$ were a trivial cofibration, this
would follow from (TM4), as in [QH, I.1.10].
In the general case, factor $w\sqcup 1$ as
$A'\sqcup A\rightarrowtail A''
\arrowsim A$.
Since $A$ and $A'$ are cofibrant, so are
the evident maps $A\rightarrowtail A''$
and $A'\rightarrowtail A''$.
They are weak equivalences, since their
composition with $A''\arrowsim A$ is.
By the trivial cofibration case, we have
the desired isomorphism.\qquad $\square$
$$
\matrix
\pi(A',B) &\longmapleft{20}{\cong} 
&\pi(A'',B) &\longmapright{20}{\cong} &\pi(A,B)\\
\noalign{\medskip}
\hskip30pt\raise6pt\hbox{$w^*$}\kern-10pt
\rlap{$\vbox{\beginpicture
\setcoordinatesystem units <.60cm,.60cm>
\linethickness=1pt
\setlinear
%
%
\linethickness= 0.500pt
\setplotsymbol ({\thinlinefont .})
\plot 11.906 18.098  9.366 20.320 /
%
%
\plot  9.599 20.201  9.366 20.320  9.516 20.105 /
\linethickness=0pt
\putrectangle corners at  9.341 20.345 and 11.932 18.072
\endpicture}$}
&&\raise15pt\hbox{$\Bigg\uparrow$} 
&&\null\kern-35pt\llap{$\vbox{\beginpicture
\setcoordinatesystem units <.55cm,.55cm>
\linethickness=1pt
\setlinear
%
%
\linethickness= 0.500pt
\setplotsymbol ({\thinlinefont .})
\plot 11.906 18.098 14.129 20.637 /
%
%
\plot 14.009 20.405 14.129 20.637 13.914 20.488 /
\linethickness=0pt
\putrectangle corners at 11.881 20.663 and 14.154 18.072
\endpicture}$}\null\kern-10pt
  \raise8pt\hbox{$\ssize =$}\\
&&\pi(A,B) &
\endmatrix
$$

\Proclaim{Lemma 5.5}
Let $\scrC$ be a basic model category, and
$A$ any object.
Suppose given weak equivalences $w'\:A'\arrowsim
A$, $w''\:A''\arrowsim A$ with $A'$ and
$A''$ cofibrant, and a fibrant object $B$.
Then for any two maps $a_1,a_2\:A'\to A''$
so that $w'=w''a_1=w''a_2$ we have
$$
a_1^*=a_2^*\: \pi(A'',B)\to
\pi(A',B).
$$
\finishproclaim

\proof
Given $f\:A''\to B$, factor $(f,w'')$ as
$A''\arrowsim P\twoheadrightarrow B\times A$.
As $B\twoheadrightarrow *$ is a fibration
and $w'$ is a weak equivalence, its
pullback $B\times A'\to B\times A$ exists
and is a weak equivalence by (TM1).
In turn, the pullback of this map along 
$P\twoheadrightarrow B\times A$ is a weak
equivalence $Q\arrowsim P$, and
$Q\twoheadrightarrow B\times A'$ is a
fibration.
We have a diagram in which the bottom
square is a pullback, so the fold map
on $A'$ factors through the dotted arrow.
$$
\matrix
A \sqcup A' &\longmapright{30}{a_1\sqcup a_2}
  &A'' &&\\
\hbox{\vbox{\beginpicture
\setcoordinatesystem units <.30cm,.30cm>
\linethickness=1pt
\setlinear
%
%
\linethickness= 0.500pt
\setplotsymbol ({\thinlinefont .})
\setdashes < 0.1270cm>
\plot 10.795 22.225 10.795 19.685 /
\setsolid
%
%
\plot 10.731 19.939 10.795 19.685 10.859 19.939 /
\setdashes < 0.1270cm>
\linethickness=0pt
\putrectangle corners at 10.770 22.250 and 10.820 19.660
\endpicture}}
 &&\null\kern-3pt\raise7pt\hbox{$\Lmapdown{\sim}$}
  &\vbox{\llap{$\beginpicture
\setcoordinatesystem units <.37cm,.37cm>
\linethickness=1pt
\setshadesymbol ({\thinlinefont .})
\setlinear
%
%
\linethickness= 0.500pt
\setplotsymbol ({\thinlinefont .})
\plot  6.032 20.955  8.572 18.733 /
%
%
\plot  8.340 18.852  8.572 18.733  8.423 18.948 /
\linethickness=0pt
\putrectangle corners at  6.007 20.980 and  8.598 18.707
\endpicture$}}\kern-9pt\raise15pt\hbox{$\ssize f$} &\\
Q &\longmapright{25}{\sim} &P
&\maptwohead{20}{} &\llap{$B$}\\
\twoheaddownarrow &&\twoheaddownarrow
  &&\\
A' &\longmapright{25}{\sim} &A &&
\endmatrix
$$
Factor the dotted arrow as $A'\sqcup
A'\rightarrowtail C\arrowsim Q$.
Now the map $A'=A'\sqcup \emptyset\to
A'\sqcup A'\to C$ is a weak equivalence,
since its composite with $C\arrowsim P$ is
$A'\arrowsim A''\arrowsim P$.
Its composite with $C\to A'$ is the
identity on $A'$, so $C\arrowsim A'$ is
also a weak equivalence.
We have shown that the map $fa_1\sqcup
fa_2$ factors as
$$
\matrix
A \sqcup A' &\longrightarrowtail &C
&\longmapright{20}{\sim} &Q
&\longmapright{20}{\sim} &P 
  &\longmapright{20}{} &B\\
&\null\kern-10pt\raise6pt\hbox{$\ssize \nabla$}\kern-10pt
\vbox{\beginpicture
\setcoordinatesystem units <.35cm,.35cm>
\linethickness=1pt
\setshadesymbol ({\thinlinefont .})
\setlinear
%
%
\linethickness= 0.500pt
\setplotsymbol ({\thinlinefont .})
\plot  6.032 20.955  8.572 18.733 /
%
%
\plot  8.340 18.852  8.572 18.733  8.423 18.948 /
\linethickness=0pt
\putrectangle corners at  6.007 20.980 and  8.598 18.707
\endpicture}
&\null\raise10pt\hbox{\,\,\,$\rmapdown{\sim}$} 
  &&&&&&\\
\vspace{-5pt}
&&A' &&&&&&
\endmatrix
$$
which shows that $fa_1\simeq_L fa_2$, as
required.\qquad $\square$

In order to have an absolutely functorial
set $[A,B]$, it is useful to parametrize
the cofibrant resolutions of $A$, and the
fibrant resolutions of $B$.

\Proclaim{Definition 5.6} \rm
Given $A$ in $\scrC$, let $I_A$ denote the category whose objects are
weak equivalences $A'\arrowsim A$ with $A'$ cofibrant.
There is a (unique) morphism from $A'\arrowsim A$ to $A''\arrowsim A$
in $I_A$ if there exists a (backwards!) arrow $A''\to A'$ so that
$A''\to A'\to A$ equals $A''\to A$.

Dually, given $B$ in $\scrC$, let $J_B$ denote the category of weak
equivalences $B\arrowsim B'$ with $B'$ fibrant, whose opposite
category is $I_{(A^{\op})}$.
\finishproclaim

\Proclaim{Lemma 5.7}
(a)\enspace
 $I_A$ is a cofiltering category with
at most one morphism between any two objects.

(b)\enspace
Given $A$ and $B$, there is an induced
functor $\pi(\,\,,\,\,)\: I_A^{\op}\times
J_B\to\Sets$.
Moreover, each morphism in
$I_A^{\op}\times J_B$ induces a bijection
on $\pi(\,\,,\,\,)$.
\finishproclaim
\goodbreak

\proof
(a) is obvious, once we note that for each
$A'$ and $A''$ we can factor $A\sqcup
A''\rightarrowtail A'''\arrowsim A$.
Part (b) follows from 5.4 and 5.5.\qquad
$\square$

\Proclaim{Definition 5.8} \rm
Given $A$ and $B$ in $\scrC$, set
$[A,B]=\lim\limits_{I_A\times
J_B}\pi(A',B')$.
\finishproclaim

\Proclaim{Lemma 5.9}
If $A$ and $B$ are cofibrant, and $C'$ is
fibrant, then a weak equivalence
$w\:B\arrowsim B'$ determines a
well-defined ``composition'' pairing
$$
\pi(A,B')\times \pi(B,C')\to\pi(A,C').
$$
\finishproclaim

\proof
Suppose given $f\:A\to B'$ and $g\:B\to
C'$.
Since $A$ and $B$ are cofibrant, we can
form $f\sqcup w\:A\sqcup B\to B'$.
Factoring it as $A\sqcup B\rightarrowtail
B''\arrowsim B'$, axiom (TM2) shows that
the composite $B\rightarrowtail A\sqcup
B\rightarrowtail B''$ is a trivial
cofibration.
Since $C'$ is fibrant, $g$ lifts to
$g''\:B''\to C'$ by (TM4).
Since $\pi(B'',C')\cong\pi(B,C')$ by 5.4,
$g''$ is well-defined up to $\simeq_L$.
We let $g\circ f$ be the composite 
of $A\rightarrowtail A\sqcup
B\rightarrowtail B''$ and $g''$.
$$
\matrix
&&B &\mapright{g} &C'\\
&&\Ldownarrowtail{\sim} 
  &\vbox{\beginpicture
\setcoordinatesystem units <.25cm,.25cm>
\linethickness=1pt
\setlinear
%
%
\linethickness= 0.500pt
\setplotsymbol ({\thinlinefont .})
\setdashes < 0.1270cm>
\plot  8.604 16.796 10.890 18.161 /
\setsolid
%
%
\plot 10.705 17.976 10.890 18.161 10.640 18.085 /
\setdashes < 0.1270cm>
\linethickness=0pt
\putrectangle corners at  8.579 18.186 and 10.916 16.770
\endpicture} &\twoheaddownarrow\\
A\sqcup B &\longrightarrowtail &B'' 
  &\longmapright{20}{} &*\\
\vspace{5pt}
\uparrowtail &&\rmapdown{\sim} &&\\
A &\longmaprightsub{20}{f} &B' &&
\endmatrix
$$
Replacing $A$ by a cylinder object in this
argument shows that if $f\simeq_L f'$ then
$g\circ f\simeq_L g\circ f'$.
Dually, if $g\simeq_R g'$ then $g\circ
f\simeq_R g'\circ f$.
\qquad $\square$
 
\Proclaim{Porism 5.10}
Suppose given $f'\:A\to B'$ and $g\:B\to
C'$, with $w\:B\arrowsim B'$.

\Item{\sans (a)}
If $f=wf_0$ for some $f_0\: A\to B$, then
$g\circ f$ is the class of $gf_0$ in
$\pi(A,C')$.

\Item{\sans (b)}
If $g=g'w$ for some $g'\:B'\to C'$, then
$g\circ f$ is the class of $g'f$ in
$\pi(A,C')$.

\noindent
In particular, $g\circ w=g$ and $w\circ
f=f$.
Hence the class $1_B$ of $w$ in $\pi(B,B)$
is a $2$-sided unit for the composition
pairing.
\finishproclaim

\Proclaim{Lemma 5.11}
The composition pairing is associative.
\finishproclaim

\proof
Suppose given $f\:A\to B'$, $g\:B\to C'$
and $h\:C\to D'$ with $A$, $B$ cofibrant and 
$B'$, $C'$, $D'$ fibrant, and weak equivalences 
$w\:B\arrowsim B'$ and $v\:C\arrowsim C'$.
The proof of 5.9 shows that we can replace
$B$ by $B''$ and $g$ by $g''$ to assume
that $f=wf_0$ for some $f_0\: A\to B$.
Similarly, we may replace $C$ by $C''$ and
$h$ by $h''$ to assume that $g=vg_0$ for
some $g_0\:B\to C$.
But then the formula
$$
(h\circ g)\circ f=(hg_0)\circ f=
(hg_0)f_0=h(g_0f_0)=h\circ(gf_0)=
h\circ (g\circ f)
$$
holds in $\pi(A,D')$, as required.\qquad
$\square$

\Proclaim{Definition 5.12} \rm
Suppose that $\scrC$ is a basic model category.
Then $\Ho\scrC$ is the category with the same objects as $\scrC$, but
the Hom-set of morphisms $A\to B$ in $\Ho\scrC$ is $[A,B]$.
Composition is defined in 5.9, and is
associative with identity by 5.10 and 5.11.

Also, there is a canonical functor $\scrC\to\Ho\scrC$, defined by~5.10.

\Proclaim{Theorem 5.13} \rm
If $\scrC$ is a basic model category then $\Ho\scrC$ is the
localization $\scrC[w^{-1}]$ of $\scrC$ at the family 
$w=\we(\scrC)$ of weak equivalences.
\finishproclaim
\proof
We first show that every weak equivalence $f\: A\arrowsim B$ in
$\scrC$ becomes an isomorphism in $\Ho\scrC$. Choosing weak
equivalences $A'\arrowsim A$ and $B\arrowsim B'$ with $A'$ cofibrant
and $B'$ fibrant, we see that $\pi(A',B')$ is canonically isomorphic
to $[A,A]$, $[A,B]$, $[B,B]$ and $[B,A]$. The class of the composite
map $f'\: A'\to B'$ represents $1_A$, $[f]$, $1_B$ and a morphism
$g\in[B,A]$, respectively. The composition pairing
$$
\pi(A',B') \times \pi(A',B') \to \pi(A',B')
$$
satisfies $f'\circ f'=f'$ by~5.10. By definition 5.8, this shows that 
$[f]\circ g=1_B$ and $g\circ [f]=1_A$, i.e., 
$[f]$ is an isomorphism in $\Ho\scrC$.

Conversely, suppose that $\phi\:\scrC\to\scrC'$ is a functor sending
weak equivalences to isomorphisms. If $C$ is a cylinder object for $A$
then the two composites $i_0,i_1\: A\to A\sqcup A \to C$
are both left inverses to $q\: C\arrowsim A$. Hence
$\phi(i_0)=\phi(i_1)$. Therefore if $f\simeq_Lg$
via a map $H\:C\to B$ (notation as in~5.1), then
$\phi(f)=\phi(g)$ because
$$
\phi(f) = \phi(H)\phi(i_0) = \phi(H)\phi(i_1) = \phi(g).
$$
It follows that $\phi$ factors uniquely through
$\scrC\to\Ho\scrC$. \qquad  
$\square$

\newpage

\centerline{\sectionfont \S6 Functor Categories}

\bigskip
If $\bfK$ is a small category and $\scrC$ is
a model category, we would like to make
the functor category $\scrC^{\bfK}$ into a
model category.
There are at least two model structures
one could consider: a local structure and
a global structure.
We shall focus on the global structure.

\Proclaim{Definition 6.1} \rm
Let $\bfK$ be a small category and $\scrC$ a
Thomason model catetory.
We say that a natural transformation
$\eta\: F\to G$ of functors $\bfK\to\scrC$ is
a ({\it global}) {\it cofibration} (written
$F\rightarrowtail G$) if for each $K$ in
$\bfK$ the map $\eta_K\:F(K)\to G(K)$ is a
cofibration in $\scrC$; {\it weak
equivalence} (written $F\arrowsim G$) if
for each $K$ in $\bfK$ the map $\eta_K\:
F(K)\to G(K)$ is a weak equivalence in
$\scrC$; ({\it global}) {\it fibration}
(written $F\twoheadrightarrow G$) for for
each $K$ in $\bfK$ the map $\eta_K$ is a
fibration in $\scrC$, and if $\eta$ has
the right lifting property with respect to
transformations 
$A\maprightarrowtail{\sim} B$ which are
simultaneously global cofibrations and weak
equivalences.
That is, every diagram of the following
form admits a factorization $B\to F$
making the resulting diagram commute.
$$
\CD A @>>> F\\
\Ldownarrowtail{\sim} @. @V{\eta}VV\\
B @>>> G
\endCD
$$
\finishproclaim

We shall say that $F\to G$ is a {\it
pointwise fibration} (resp. has pointwise
$P$) if each $F(i)\to G(i)$ is a fibration
(resp. has property $P$) in $\scr C$.

\Proclaim{Proposition 6.2 {\rm [T85, 32]}}
If $\scrC$ is a basic model category and
$\bfK$ is a small category, then $\scrC^{\bfK}$
with the global structure above satisfies
axioms (TM0)--(TM3).
If $\scrC$ satisfies the factorization
axiom (TM5c) then $\scrC^{\bfK}$ satisfies
(TM4) and (TM5c).
\finishproclaim

\proof
The constant functors $F(K)=\emptyset$ and
$F(K)=*$ are the initial and terminal
objects of $\scrC^{\bfK}$.
Axioms (TM0), (TM2) and (TM3) for
$\scrC^{\bfK}$ are immediate from the pointwise
definitions.
Because limits and colimits in $\scrC^{\bfK}$
are defined pointwise, there is only one
nontrivial point to check for (TM1): if
$F\rightarrowtail G$ is a fibration in
$\scrC^{\bfK}$, does the pullback
$\eta\:F\times_G H\to H$ have the lifting
property with respect to a trivial
cofibration $A\maprightarrowtail{\sim}B$?
To see that it does, consider the diagram
$$
\CD
A @>>> F\times_G H @>>> F\\
\Ldownarrowtail{\sim} @. @VVV
  &\Rtwoheadmapdown{\eta}\\
B @>>> H @>>> G.
\endCD
$$
The lifting property for the global
fibration $\eta$ yields a fill-in $B\to F$
which must factor through a map to the
pullback $B\to F\times_G H$.
So (TM1) holds for $\scrC^{\bfK}$.

If $\scrC$ satisfies (TM5c), any $\eta\:
A\to B$ in $\scrC^{\bfK}$ determines
factorizations which are functorial in
$K$:
$$
A(K)\longrightarrowtail T(\eta_K)\arrowsim 
B(K).
$$
This gives the functorial factorization
$A\longrightarrowtail T\arrowsim B$ in
$\scrC^{\bfK}$.

To check axiom (TM4), consider a square in
$\scrC^I$ of the form
$$
\CD
A @>>> X\\
\downarrowtail @. \twoheaddownarrow\\
B @>>> Y.
\endCD
$$
If $A\arrowsim B$ is a weak equivalence,
the definition of global fibration yields
a fill-in $B\to X$.
Suppose then that $X\arrowsim Y$ is a weak
equivalence.
By (TM1) the pushout $X\rightarrowtail C$
of $A\rightarrowtail B$ exists, and there
is a canonical map $C\to Y$.
By (LM5) we can factor it as
$C\rightarrowtail D\arrowsim Y$.
By saturation (TM2), the cofibration
$X\rightarrowtail D$ is also a weak
equivalence since $D\arrowsim Y$ and
$X\arrowsim Y$ are.
But then the definition of global
fibration yields a fill-in $D\to X$ of the
right square in the diagram
$$
\CD
A @>>> X @= X\\
\downarrowtail @. \Ldownarrowtail{\sim}
 @. \Rtwoheadmapdown{\sim}\\
B @>>> D @>{\sim}>> Y,
\endCD
$$
and $B\to D\to X$ provides the desired
factorization.
This completes the verification of
(TM4).\qquad $\square$

\Proclaim{Remark 6.3} \rm
If $\scrC$ is a Thomason model category,
all that is missing is the factorization
part (f) of (CM5).
We do have a functorial factorization
$F\arrowsim M\longrightarrow
 G$ from (RM5) on $\scrC$, and $M\to G$ is
a pointwise fibration, i.e., each
$M(K)\twoheadrightarrow G(K)$ is a
fibration in $\scrC$.
However, there appears to be no reason for
$M\to G$ to satisfy the right lifting
property required of a global fibration.

Recall that an object $C$ of $\scrC$ is called {\it fibrant} if
$C\twoheadrightarrow *$ is a fibration.
We write $\scrC_{\fib}$ for the full
subcategory of fibrant objects in $\scrC$
and $\scrC_{\fib}^{\bfK}$ for the category of
all pointwise fibrant functors, i.e., all
functors ${\bfK}\to\scrC_{\fib}$.
Note that $(\scrC^{\bfK})_{\fib}$ is strictly
contained in
$\scrC_{\fib}^{\bfK}=(\scrC_{\fib})^{\bfK}$.
The following argument allows us to reduce
the verification of the fibration half 
(TM5f) of axiom (TM5) to pointwise
fibrations in $\scrC_{\fib}^{\bfK}$.
\finishproclaim

\Proclaim{Lemma 6.4 {\rm [T85, 166]}}
Suppose that there is a functorial
factorization in $\scrC^{\bfK}$
$$
F'\arrowsim H'\twoheadrightarrow G'
$$
for every pointwise fibration $F'\to G'$
of pointwise fibrant functors.
Then axiom (TM5f) holds in $\scrC^{\bfK}$, and
$\scrC^{\bfK}$ is a Thomason model category.
\finishproclaim

\proof
Suppose given $F\to G$ in $\scrC^{\bfK}$.
Using functorial factorization in $\scrC$,
we have a factorization $G\arrowsim
G'\longrightarrow *$, where $G'$ is
pointwise fibrant, and then a
factorization of $F\to G'$
$$
\CD
F @>>> G\\
\lmapdown{\sim} @. \lmapdown{\sim}\\
F' @>>> G'
\endCD
$$
in which $F'\to G'$ is pointwise fibrant.
By assumption, it factors as $F'\arrowsim
H'\twoheadrightarrow G'$.
By (TM1) the pullback $H=F'\times_{G'} G$
exists, $H\twoheadrightarrow G$ is a
fibration and $H\arrowsim H'$ is a weak
equivalence.
By saturation (TM2), $F\arrowsim H$ is a
weak equivalence, and we have the required
functorial factorization $F\arrowsim
H\twoheadrightarrow G$.\qquad $\square$
$$
\CD
F @>>> H @. \maptwohead{15}{} &G\\
@V{\sim}VV @VV{\sim}V  &@V{\sim}VV\\
F' @>{\sim}>> H'  @.\maptwohead{5}{} &G'
\endCD
$$

\newpage

\centerline{\sectionfont \S7 Enriched
Functor Categories}

\bigskip
If $\scrC$ is a Thomason model category,
we saw in the last section that
$\scrC^{\bfK}=\Cat(\bfK,\scrC)$ satisfies
all the axioms for a model category except
for (TM5f), that every map $F\to G$ in $\scrC^{\bfK}$ has a
factorization $F\arrowsim H\twoheadrightarrow G$.
In order to do this we shall use homotopy
limits and homotopy ends in $\scrC$; this requires that
$\scrC$ be not only complete but enriched
over finite simplicial sets.
Recall from 1.12 that this means that
$\scrC$ has a mapping object functor
satisfying axioms (RE1) -- (RE5).
We will also need $\scrC$ to satisfy the
following fibration conditions.

\Proclaim{Definition 7.1 {\rm [T85, 147]}} \rm
We say that {\it towers preserve
fibrations} in $\scrC$ if given a map
between towers of fibrations in $\scrC$
$$
\matrix
\cdots &\twoheadrightarrow  &C_n
&\twoheadrightarrow &C_{n-1}
&\twoheadrightarrow &\cdots
&\twoheadrightarrow &C_0
&\twoheadrightarrow &*\\
&&\Big\downarrow &&\Big\downarrow
  &&&&\twoheaddownarrow &&\\
\cdots &\twoheadrightarrow &D_n 
&\twoheadrightarrow &D_{n-1} 
&\twoheadrightarrow &\cdots
&\twoheadrightarrow &D_0
&\twoheadrightarrow &*
\endmatrix
$$
in which $C_0\to D_0$ is a fibration and each
$C_{n+1}\twoheadrightarrow C_n\times_{D_n}D_{n+1}$ is a fibration 
($n\Ge 0$), then $\limleft\,C_n\twoheadrightarrow
\limleft\,D_n$ is a fibration.
If in addition each $C_n\arrowsim D_n$ is a
weak equivalence, so is
$\limleft\,C_n\maptwohead{5}{\sim}
\limleft\,D_n$.
Note that this implies that
$\limleft\,C_n$ is fibrant (take $D_n=*$).

We say that {\it products preserve
fibrations} ({\it and equifibrations}) in
$\scrC$ if for every set $I$ and every
family $\{C_i\}_{i\in I}$ of fibrant
objects of $\scrC$: (a) $\prod C_i$ exists
in $\scrC$, and (b) if
$B_i\twoheadrightarrow C_i$ are fibrations
(resp. equifibrations) then $\prod
B_i\twoheadrightarrow \prod C_i$ is a
fibration (resp. an equifibration) in
$\scrC$.
Note that (b) applied to
$C_i\twoheadrightarrow *$ implies that
$\prod C_i$ is fibrant.
\finishproclaim
\goodbreak

\Proclaim{Theorem 7.2 {\rm [T85, 160]}}
Let $\scrC$ be a complete Thomason model
category with a right enrichment over
finite simplicial sets (1.12).
Suppose that products and towers in
$\scrC$ preserve fibrations and
equifibrations.
If $\bfK$ is a small category then
$\scrC^{\bfK}$ has the structure of a
Thomason model category.

Moreover, $\scrC^{\bfK}$ is right
enriched, and products and towers in
$\scrC^{\bfK}$ preserve fibrations and
equifibrations.
\finishproclaim

The proof of theorem 7.2 will be
 given after lemma 7.6 below.
It requires a simplicial $hom$ construction. 

\Proclaim{Definition 7.3} \rm
Given $B$ and $C$ in $\scrC$, the simplicial set
$\hom_{\scrC}(B,C)$ is defined to be
$[\hom_{\scrC}(B,C)]_p=\Hom_{\scrC}(B,\Map(\Delta
[p],C))$.
Thus $\hom_{\scrC}$ is a functor from
$\scrC^{\op}\times\scrC$ to
simplicial sets, and we may take its
homotopy end.
\finishproclaim

\Proclaim{Lemma 7.4 {\rm [T85, 121]}}
If $A\rightarrowtail B$ is a cofibration
in $\scrC$, then
$$
\hom_{\scrC}(B,C)\twoheadrightarrow
\hom_{\scrC} (A,C)
$$
is a Kan fibration of simplicial sets.
\finishproclaim

\proof
For $0\Le k\le n$, let $V=V(n,k)$ denote the
simplicial subcomplex of $\Delta[n]$
which is the union of all faces except the
$k$th one.
By axiom (RE2) the right vertical map is a
fibration in $\scrC$ in the diagram:
$$
\matrix
A &\longrightarrow &\Map(\Delta[n],C)\hfill\\
\downarrowtail &\vbox{\beginpicture
\setcoordinatesystem units <.30cm,.30cm>
\linethickness=1pt
\setlinear
%
%
\linethickness= 0.500pt
\setplotsymbol ({\thinlinefont .})
\setdashes < 0.1270cm>
\plot  8.604 16.796 10.890 18.161 /
\setsolid
%
%
\plot 10.705 17.976 10.890 18.161 10.640 18.085 /
\setdashes < 0.1270cm>
\linethickness=0pt
\putrectangle corners at  8.579 18.186 and 10.916 16.770
\endpicture} &\twoheaddownarrow\\
B &\longrightarrow &\Map(V(n,k),C).
\endmatrix
$$
By (TM4), a fill-in exists.
Thus the function
$$
\align
&\Hom_{\scrC}(B,\Map(\Delta[n],C))\to\\
&\qquad\qquad\qquad
 \Hom_{\scrC}(A,\Map(\Delta[n],C))
\operatornamewithlimits{\times}
\limits_{\Hom_{\scrC}(B,\Map(V,C))}
 \Hom_{\scrC}(A,\Map(V,C))
\endalign
$$
is surjective.
By (RE2), $\Map(-,C)$ sends the pushout
$V(m,k)$ of the cofibrations
$\Delta[n-2]\rightarrowtail \Delta[n-1]$ to
the pullback in $\scrC$ of the fibrations
$\Map(\Delta[n-1],C)\twoheadrightarrow
\Map(\Delta[n-2],C)$.
Applying $\Hom_{\scrC}(B,-)$, we see that
$\Hom_{\scrC}(B,\Map(V,C))$ is the
set-theoretic pullback of functions
$\hom(B,C)_{n-1}\to\hom(B,C)_{n-2}$.
That is,
$$
\Hom_{\scrC}(B,\Map(V,C))\cong
\Hom_{\Delta^{\op}\Sets}(V,\hom_{\scrC}(B,C)).
$$
Since
$\Hom_{\scrC}(B,\Map(\Delta[n],C))=\hom_{\scrC}
(B,C)_n$ by definition, 
this shows that\break
 $\hom_{\scrC}(B,C)
\to \hom_{\scrC}(A,C)$ is a Kan fibration.
\qquad $\square$.

If $C$ is fibrant in $\scrC$, then each
$\Map(\Delta[p],C)$ is fibrant by (RE1).

\Proclaim{Theorem 7.5 {\rm [T85, 158]}}
For $B$ in $\scrC$ and
$F\:\bfK\to\scrC_{\fib}$, there is a
natural bijection:
$$
\Hom_{\scrC}\left(
B,\operatornamewithlimits{\holim}
  \limits_{\bfK} F\right)\cong
\left[\operatornamewithlimits{\holim}\limits
  _{K\in\bfK}\hom_{\scrC}(B,FK)\right]_0.
$$
For $B$ in $\scrC^{\bfK}$, there is a
natural bijection of sets:
$$
\Hom_{\scrC^{\bfK}}\left(B, \holim\,F \right)
\cong \left[\hoend
\hom_{\scrC}(BK',FK)\right]_0.
$$
\finishproclaim 

\proofof\hbox{([T85, 123--130]):}
By 1.4 and 1.14, part a) is a special case of b).
By [Mac, p.~219], the left side is the end
$$
\int\nolimits_{\bfK}\Hom_{\scrC}\left(BK,
\holim\limits_{K/\bfK}\,F\right)=
\int\nolimits_{\bfK}\Hom_{\scrC}
\left(BK,\int\nolimits_{\Delta}\Map
\Bigl(\Delta[p],\prod\limits_{I}FK_p\Bigr)
\right),
$$
where the indexing set $I$ runs over all
diagrams $K\to K_0\to\cdots\to K_p$ in
$\bfK$.
By the universal property of ends, this
equals the double end
$$
\gather
\int\limits_{\bfK}\int\limits_{\Delta}
\Hom_{\scrC}\left(BK,\Map\Bigl(\Delta[p],
\prod\limits_{I} FK_p\Bigr)\right).\\
\intertext{By the definition of
$\hom_{\scrC}$, this equals}
\iint\limits_{\bfK\times\Delta}
\hom_{\scrC}\Bigl(BK,\prod\limits_{I}
FK_p\Bigr)_p
=\iint\limits_{\bfK\times\Delta}
\Hom_{\Delta^{\op}\Sets}\left(\Delta[p],
\hom_{\scrC}\Bigl(BK,\prod\limits_{I}FK_p\Bigr)
\right).
\endgather
$$
Because $\Hom_{\Delta^{\op}\Sets}(X,Y)$ is
the set of $0$-simplices in $\Map(X,Y)$,
and limits in $\Delta^{\op}\Sets$ are
formed pointwise, this equals
$$
\iint\limits_{\bfK\times\Delta}
[\Map(\Delta[p],Y)]_0=
\biggl[\,\,\,\,\iint\limits_{\bfK\times\Delta}
\Map(\Delta[p],Y)\biggr]_0,
$$
where $Y=Y (K',K)$ denotes
$\hom_{\scrC}
\left(BK,\prod\nolimits_{I}FK_p\right)$.
Hence it suffices to show that we have an
isomorphism of simplicial sets between
$$
\iint\limits_{\bfK\times\Delta}\Map(\Delta
[p],Y)\qquad\text{and}\qquad
\hoend\hom_{\scrC}(BK',FK).
$$
Because $\Map(\Delta[p],-)$ preserves
limits of simplicial sets and ends
commute, the left side equals
$$
\int\limits_{\Delta}\int\limits_{\bfK}
\Map(\Delta[p],Y)=\int\limits_{\Delta}
\Map\left(\Delta[p],\int\nolimits_{\bfK}Y\right).
$$
Now axiom (RE5) in $\scrC$ implies that
(for each $p$ and $K$)
$$
Y(K',K)=\hom_{\scrC}\biggl(BK',
\prod\limits_{K\to\cdots\to
K_p}FK_p\biggl)\cong\prod\limits_{K\to\cdots\to
K_p}\hom_{\scrC}(BK',FK_p).
$$
But now $\hoend\hom_{\scrC}(BK',FK)$ is
defined (1.3) to be the simplicial set
$$
\int\nolimits_{\Delta}\Map\left(\Delta[p],
\prod\limits_{K_0\to\cdots\to
K_p}\hom_{\scrC}(BK_0,FK_p)\right).
$$
Hence it suffices to show that for all $p$
$$
\prod\limits_{K_0\to\cdots\to K_p}
\hom_{\scrC}(BK_0,FK_p)\cong
\int\limits_{\bfK}\prod\limits_{K\to
K_0\to \cdots\to K_p}\hom_{\scrC}(BK',FK_p).
$$

We can apply lemma 7.6 below to
$I=\bfK^{\op}$, $i_0=K$, $\scrS=$
simplicial sets and
$G(K')=\prod\limits_{K_0\to\cdots\to
K_p}\hom_{\scrC}(BK',FK_p)$, where the
indexing set runs over diagrams of length
$p$ in $\bfK$.
The result is the required isomorphism:
$$
\alignat2
G(K_0) &\cong \int\nolimits_{\bfK^{\op}}
  (K,K')\longmapsto\prod\limits
\Sb
K\to K_0\\
\text{in }\bfK
\endSb
G(K') &\qquad &\\
&\cong\int\nolimits_{\bfK}(K',K)\longmapsto
  \prod\limits_{K\to K_0}G(K') &\qquad &\\
&=\int\nolimits_{\bfK}(K',K)\longmapsto
\prod\limits_{K\to K_0
  \to\cdots\to K_p}\hom_{\scrC}(BK',
  FK_0). &\qquad &\square
\endalignat
$$

\Proclaim{Lemma 7.6 {\rm [T85, 129]}}
Let $I$ be a small category, $i_0$ an
object of $I$, and $G\: I\to\scrS$ a
functor to a complete category $\scrS$.
If $F\: I^{\op}\times I\to\scrS$ denotes
the functor $F(i,j)=\prod\limits_{i_0\to
i}G(j)$ then there is a natural
isomorphism
$$
G(i_0)\cong\int\nolimits_{I} F=
\int\nolimits_{I}(i,j)\longmapsto
\prod\limits_{i_0\to i}G(j).
$$
\finishproclaim

\proof
The diagonal
$G(i_0)\to\prod\limits_{i_0\to i}G(i_0)\to
\prod\limits_{i_0\to i}G(i)$ is dinatural,
and induces a map $G(i_0)\to
\int\nolimits_{I}F$.
Conversely, the projection onto the
coordinate indexed by the identity of
$i_0$ yields a map $\int\nolimits_{I}
F\to\prod\limits_{i_0\to i_0}G(i_0)\to
G(i_0)$.
These maps are inverses to each other.
\qquad $\square$

\proofof\hbox{(of theorem 7.2):}
As remarked above, we are reduced by 6.2
and 6.4 to producing a functorial
factorization of $\eta\:F\to G$ under the
assumption that each $FK\twoheadrightarrow
GK\twoheadrightarrow *$ is a fibration.
 From the diagram using lemma 1.13:
$$
\CD
F @>{\sim}>> \holim F\\
@V{\eta}VV @VV{\holim \eta}V\\
G @>{\sim}>> \holim G.
\endCD
$$

Now $\holim \eta$ is a pointwise
fibration by lemma 1.8.
Let $H$ denote the pullback of $G$ along
$\holim \eta$; axiom (TM1) in $\scrC$
implies that $H\arrowsim \holim F$ is a
weak equivalence and $H\to G$ is a
pointwise fibration.
Saturation (TM2) gives $F\arrowsim H$.
We claim that in the resulting functorial
factorization of $\eta$,
$$
F\arrowsim H\longrightarrow G,
$$
the map $H\to G$ is a global fibration,
i.e., it satisfies the right lifting
property.

Let $A\maprightarrowtail{\sim} B$ be a
cofibration and weak equivalence.
We need to show that a fill-in exists for
any diagram
$$
\CD
A @>>> H @>>> \holim F\\
\Ldownarrowtail{\sim} &&@VVV @VV{\holim\eta}V\\
B @>>> G @>>> \holim G.
\endCD
$$
As the right square is a pullback, it
suffices to find a fill-in $B\to\holim F$,
i.e., to show that $\holim\eta$ is a
global fibration.
That is, it suffices to show that the
following map is a surjection.
$$
\Hom_{\scrC^{\bfK}}(B,\holim F)\to
\Hom_{\scrC^{\bfK}}(A,\holim
F)]\operatornamewithlimits{\times}
  \limits_{\Hom(A,\holim G)}
\Hom_{\scrC^{\bfK}}
(B,\holim G).
$$

By proposition 7.5, $\Hom(B,\holim F)$
is naturally isomorphic to the degree zero
part of the homotopy end 
$\hoend_{\bfK} \hom_{\scrC}(BK',FK)$,
because each $FK$ is fibrant.

As the homotopy end formed in
$\Delta^{\op}$-$\Sets$ preserves limits,
such as pullbacks along a fibration, the
map we need to show surjective is
identified with the degree zero part of the
simplicial map
$$
\align
\hoend &\hom_{\scrC}(BK',FK)\\ 
&\to \left[\hoend \hom(AK',FK)\right]
  \times_{\hoend\hom(AK',GK)}
\left[\hoend\hom(BK',GK)\right]\\
&\cong
\hoend\Bigl[\hom(AK',FK)\times_{\hom(AK',GK)}
\hom(BK',GK)\Bigr].
\endalign
$$
Because each $AK'\rightarrowtail BK'$ is a
cofibration, each simplicial map
$$
\hom_{\scrC}(BK',GK)\to\hom_{\scrC}(AK',GK)
$$
is a Kan fibration by lemma 7.4.
Since each $FK\twoheadrightarrow
GK\twoheadrightarrow *$ is a fibration,
the map
$$
\hom_{\scrC}(BK',FK)\to\hom(AK',FK)
\times_{\hom(AK',GK)}\hom(BK',GK)
$$
is a fibration and a weak equivalence of
simplicial sets.
By lemma 1.8, applying $\hoend$
yields a Kan fibration and weak equivalence
of simplicial sets.
But any such map is surjective on
$0$-simplices.
This completes the proof that (TM5f) holds,
and that $\scrC^{\bfK}$ is a Thomason
model category.

To see that products preserve fibrations and
equifibrations in $\scrC^{\bfK}$,
suppose given a family  $\{F_i\twoheadrightarrow
G_i\}_{i\in I}$ of fibrations (resp.
equifibrations).
Since $\scrC$ is complete, the maps $\prod
F_i K\twoheadrightarrow \prod G_i K$ exist and are
fibrations (resp. equifibrations) in
$\scrC$ for all $K$.
It suffices to show that $\prod F_i\to
\prod G_i$ has the right lifting property.
But for any $A\maprightarrowtail{\sim}B$,
each square
$$
\CD
A @>>> \prod F_i\\
\Ldownarrowtail{} @. @VVV\\
B @>>> \prod G_i
\endCD
$$
will admit a lift since lifts $B\to F_i$
exist for each $i\in I$ and $\prod F_i$ is
a product.

To see that towers preserve fibrations in
$\Cat(\bfK,\scrC)$, suppose given a map
between towers of fibrations in which
$F_0\twoheadrightarrow G_0$ and all
$F_{n+1}\twoheadrightarrow
G_{n+1}\times_{G_n}F_n$ are fibrations.
$$
\matrix
\cdots &\twoheadrightarrow &F_n
&\twoheadrightarrow &F_{n-1}
&\twoheadrightarrow &\cdots
&\twoheadrightarrow &F_0
&\twoheadrightarrow &*\\
&&\Big\downarrow &&\Big\downarrow
&&&&\Big\downarrow &&\\
\cdots &\twoheadrightarrow &G_n 
&\twoheadrightarrow &G_{n-1} 
&\twoheadrightarrow &\cdots
&\twoheadrightarrow &G_0
&\twoheadrightarrow &*
\endmatrix
$$
Since towers preserve fibrations in
$\scrC$, each $\limleft F_nK$ and
$\limleft G_nK$ is fibrant and $\limleft
F_nK\to \limleft G_nK$ is a fibration (and
is a weak equivalence when each $F_nK\to
G_nK$ is).
To show it is a fibration it suffices to
check the lifting property.
Let $A\maprightarrowtail{\sim}B$ be an
equi-cofibration and suppose given a
square
$$
\CD
A @>>> \limleft F_n @>>> F_n\\
\Ldownarrowtail{\sim} &&@VVV @VVV\\
B @>>> \limleft G_n @>>> G_n
\endCD
$$
We need to find a lift $B\to\limleft F_n$.

We proceed by induction to produce a
coherent tower of lifts $B\to F_n$;
passing to the inverse limit will yield
the desired lift $B\to\limleft F_n$.
For $n=0$, a lift $B\to F_0$ exists
because $F_0\to G_0$ is a fibration.
Inductively, there is a lift $B\to F_n$.
It factors uniquely through the pullback
$G_{n+1}\times_{G_n}F_n$, and we have a
diagram
$$
\matrix
A &\longrightarrow 
&F_{n+1} &\longrightarrow &F_n\\
\Ldownarrowtail{} &&\twoheaddownarrow 
&\vbox{\beginpicture
\setcoordinatesystem units <.25cm,.25cm>
\linethickness=1pt
\setlinear
%
%
\linethickness= 0.500pt
\setplotsymbol ({\thinlinefont .})
\plot 11.906 18.098 14.863 19.723 /
%
%
\plot 14.671 19.545 14.863 19.723 14.610 19.656 /
\linethickness=0pt
\putrectangle corners at 11.881 19.748 and 14.889 18.072
\endpicture} &@VVV\\
B &\longrightarrow
 &G_{n+1}\times_{G_n}F_n &\longrightarrow &G_n.
\endmatrix
$$
Since the middle vertical is a fibration,
there is a lift $B\to F_{n+1}$ compatible
with the preceeding lift $B\to F_n$.
This completes the inductive step, showing
that the map $\limleft F_n\to \limleft
G_n$ is a fibration, and that towers
preserve fibrations in $\Cat(\bfK,\scrC)$.

Finally, we need to show that
$\Cat(\bfK,\scrC)$ is right enriched over $\scrS_f$.
Because $\Cat(\bfK,\scrC)_{\fib}\subseteq
\Cat(\bfK,\scrC_{\fib})$, there is a
functor
$$
\Map(\,\,,\,\,)\:\scrS_f^{\op}\times
\Cat(\bfK,\scrC)_{\fib}\longrightarrow
\Cat(\bfK,\scrC_{\fib})
$$
We need to show that this lands in
$\Cat(\bfK,\scrC)_{\fib}$ and makes
$\Cat(\bfK,\scrC)$ right enriched over
$\scrS_f$.
This is straightforward, and left to the
reader.
See [T86, 61].

This finishes the proof of Theorem 7.2
\qquad $\square$

\bigskip\bigskip
\centerline{\sectionfont Acknowledgements}

\bigskip
This work was done at the Institute for
Advanced Study in Princeton, but the main
sources were located at the Universit\'e de Paris 7. 
The author is deeply grateful to 
Liliane Barenghi, Bruno Kahn and Max Karoubi for
their help in giving me access to Thomason's private notebooks
[T78]--[T86].
I am also grateful to Luca Barbieri-Viale, Marco Grandis
and Ieke Moerdijk for providing me with their notes of seminar talks
in which Thomason presented his ideas in Genova and Utrecht.
Finally, I am grateful to Rick Jardine and Bill Dwyer for their
comments on an earlier draft.

\newpage
\references{40}

{\nspace{
\Ref[Baues]
H. Baues, 
\it Algebraic Homotopy,
\rm Cambridge Univ. Press, 1989.

\Ref[BF]
A. Bousfield and E. Friedlander,
\it Homotopy theory of $\Gamma$-spaces,
spectra and bisimplicial sets,
\rm Lecture Notes in Math. 658, 1978,
80--130.

\Ref[BK]
A. Bousfield and D. Kan,
\it Homotopy Limits, Completions and Localizations,
\rm Lecture Notes in Math. 304,
Springer-Verlag, 1972.

\Ref[Br]
K. S. Brown
\it Abstract Homotopy Theory and Generalized Sheaf Cohomology,
\rm Trans. AMS 186 (1973), 419--458.

\Ref[BW]
H. Baues and G. Wirsching,
\it The cohomology of small categories,
\rm J. Pure Appl. Alg. 38 (1985), 187--211.

\Ref[DK]
W. Dwyer and D. Kan,
\it Hochschild-Mitchell cohomology of simplicial categories
   and the cohomology of simplicial diagrams of simplicial sets,
\rm Nederl. Akad. Wetensch. Indag. Math. 50 (1988), 111--120.

\Ref[GJ]
P. Goerss and J. F. Jardine,
\it Simplicial Homotopy Theory,
\rm Birkh\"auser, 1999.

\Ref[JP]
M. Jibladze and T. Pirashvili,
\it Cohomology of algebraic theories,
\rm J. Alg. 137 (1991), 253--296.

\Ref[Mac]
S. Mac Lane,
\it Categories for the Working Mathematician,
\rm Springer-Verlag, 1971.

\Ref[Mau]
L. Mauri,
\it Thomason models for homotopy,
\rm preprint, 2000.

\Ref[PW]
T. Pirashvili and F. Waldhausen,
\it MacLane homology and topological Hoch-\break
schild homology, 
\rm J. Pure Appl. Alg. 82 (1992), 81--98.

\Ref[QH]
D. Quillen,
\it Homotopical Algebra,
\rm Lecture Notes in Math. 43,
Springer-Verlag, 1967.

\Ref[QR]
D. Quillen, 
\it Rational homotopy theory,
\rm Annals Math. 90 (1969), 205--295.

\Ref[T78]
R. Thomason, private notebook, Paris, 
Jan.-Feb., 1995, 180 pp.

\Ref[T79]
R. Thomason, private notebook, Paris,
March, 1995, 180 pp.

\Ref[T81]
R. Thomason, private notebook, Paris,
June, 1995, 180 pp.

\Ref[T82]
R. Thomason, private notebook, Paris, July
1995, 180 pp.

\Ref[T83]
R. Thomason, private notebook, Paris, 
Aug., 1995, 180 pp.

\Ref[T85]
R. Thomason, private notebook, Paris,
Sept. 1995, 180 pp.

\Ref[T86]
R. Thomason, private notebook, Paris,
Oct., 1995, 74~pp.

\Ref[TT]
R. Thomason and T. Trobaugh,
\it Higher algebraic $K$-theory of schemes
and of derived categories,
\rm The Grothendieck Festschrift III,
Progress in Math. 88, Birkh\"auser, 1990,
247--435.

\Ref[Wa]
F. Waldhausen, 
\it Algebraic $K$-theory of spaces,
\rm Lecture Notes in Math. 1126,
Springer-Verlag, 1985, 318--419.

\Ref[WT]
C. Weibel,
\it The Mathematical Enterprises of Robert
Thomason,
\rm Bull. AMS 34 (1997), 1--13.

}}

\enddocument